\documentclass{amsart}
\def\sectiontheorems{} 
\usepackage[hmargin=2cm,vmargin=2cm]{geometry} 
\usepackage{tikz-cd}
\usepackage{mdframed}
\usepackage{amssymb} 
\usepackage{biblatex}
\usepackage{hyperref}
\hypersetup{
    colorlinks=true,
    linkcolor=black,
    citecolor=black,
    urlcolor=black
}
\usepackage{mathtools}
\usepackage{mathrsfs}
\usepackage[shortlabels]{enumitem}
\usepackage[nameinlink]{cleveref}
\usepackage{comment}
\usepackage[english]{babel}
\usepackage{amsthm}
\usepackage{graphicx}
\usepackage{subcaption}
\usepackage{float}
\usepackage{empheq}
\usepackage[dvipsnames]{xcolor}
\usepackage{algorithm}
\usepackage{algpseudocode}

\theoremstyle{plain}
\ifdefined\sectiontheorems
    \newtheorem{theorem}{Theorem}[section]
\else
    \newtheorem{theorem}{Theorem}
\fi
\newtheorem{proposition}[theorem]{Proposition}
\newtheorem{corollary}[theorem]{Corollary}
\newtheorem{lemma}[theorem]{Lemma}

\newtheorem*{theorem*}{Theorem}
\theoremstyle{definition}
\newtheorem{definition}[theorem]{Definition}
\newtheorem{example}[theorem]{Example}
\newtheorem*{remark}{Remark}
\newtheorem{assumption}[theorem]{Assumption}

\renewcommand{\Im}{\operatorname{Im}}

\newcommand{\tx}{\text}

\newcommand{\dup}{\mathrm{d}} 

\title{Tracial Joint Spectral Measures For Nearly Commuting Matrices}
\author{Chenning Xu}
\address{Division of Physics, Mathematics and Astronomy, California Institute of Technology, Pasadena, CA 91125}
\email{cxu7@caltech.edu}
\date{July 2026}

\begin{document}

\maketitle

\begin{abstract}
Tracial joint spectral measures are measures on the plane and originally introduced by Hein{\"a}vaara to prove certain trace inequalities for two matrices. Hein{\"a}vaara gave complicated formulas for the density and singular part of the measure with respect to Lebesgue measure, involving the eigenvalues of an intermediate matrix. We study the density when the matrices nearly commute. We find a low dimensional approximation of the density and prove geometric properties of the support of density. Specifically, we prove the support is located near certain lines up to the first order. Up to order 3/2, the support resembles a union of ellipses, and we find formulas for the ellipses.
\end{abstract}

\section{Introduction}

Tracial joint spectral measures, introduced by Hein{\"a}vaara in 2025 \cite{Hei_spectral}, are positive measures on $\mathbb R^2$ that capture the joint spectrum of two Hermitian matrices.

\begin{theorem}[Tracial Joint Spectral Measure, {\cite[Theorem 1.4]{Hei_spectral}}]\label{thm:tracial_joint_spectral_measure_original}
For any $n\times n$ Hermitian matrices $A$ and $B$, there exists a positive measure $\mu_{A,B}$ on $\mathbb R^2$ s.t. for any positive integer $k$ and $x, y\in \mathbb R$, we have,
\begin{align*}
    \mathrm{tr}((xA+yB)^k) = k(k+1) \int_{\mathbb R^2} (ax+by)^k \mathrm{d}\mu_{A,B}.
\end{align*}
where we call $\mu_{A,B}$ the ``tracial joint spectral measure''.
\end{theorem}

As a corollary of this theorem, for any $A$ and $B$ in the Schatten-$p$ space, their span can be isometrically embedded into $L_p(\mu_{A,B})$. Additionally, Hein{\"a}vaara also found explicit formula for $\mu_{A,B}$:

\begin{theorem}[Density of Tracial Joint Spectral Measure, {\cite[Theorem 1.7]{Hei_spectral}}]
\label{thm:density_of_measure_original}
Let $n, A, B$ and $\mu_{A,B}$ be as in Theorem~\ref{thm:tracial_joint_spectral_measure_original}. Decompose $\mu_{A,B}$ into absolutely continuous part and singular part with respect to the Lebesgue measure and denote the two parts by $\mu_c$ and $\mu_s$, respectively. Assume some linear combination of $A$ and $B$ is invertible. Then the density of $\mu_c$ is given by:
\begin{align}
\label{eq:density_muc_in_intro}
    \frac{\dup\mu_c}{\dup m_2}(a,b) = \frac{1}{2\pi} \sum_{i=1}^n \left|\Im\left( \lambda_i\left(\left(I-\frac{A}{a}\right)(bA-aB)^{-1}\right)\right)\right|.
\end{align}
Furthermore, if $A$ is invertible and $A^{-1}B$ has $n$ distinct eigenvalues, the singular part $\mu_s$ satisfies 
\begin{align*}
    \mu_s(\varphi) = \sum_{v \in E(A^{-1}B)} \int_0^1 \frac{1-t}{t} \varphi(\langle Av, v\rangle t, \langle Bv, v\rangle) \dup t,
\end{align*}
where $E(C)$ denotes a set of normalized eigenvalues of an $n\times n$ complex matrix $C$ and $\varphi$ is a smooth function with compact support that does not contain 0.
\end{theorem}

We note that the formula (\ref{eq:density_muc_in_intro}) is not the original formulation in Hein{\"a}vaara’s Theorem 1.7; it is a modification based on Remark 3.2 of the paper.\\

We pick some sample matrices and plot their respective $\mu_c$ and $\mu_s$:

\begin{figure}[H]
    \centering
    \begin{subfigure}[t]{0.48\textwidth}
        \centering
        \includegraphics[width=\linewidth]{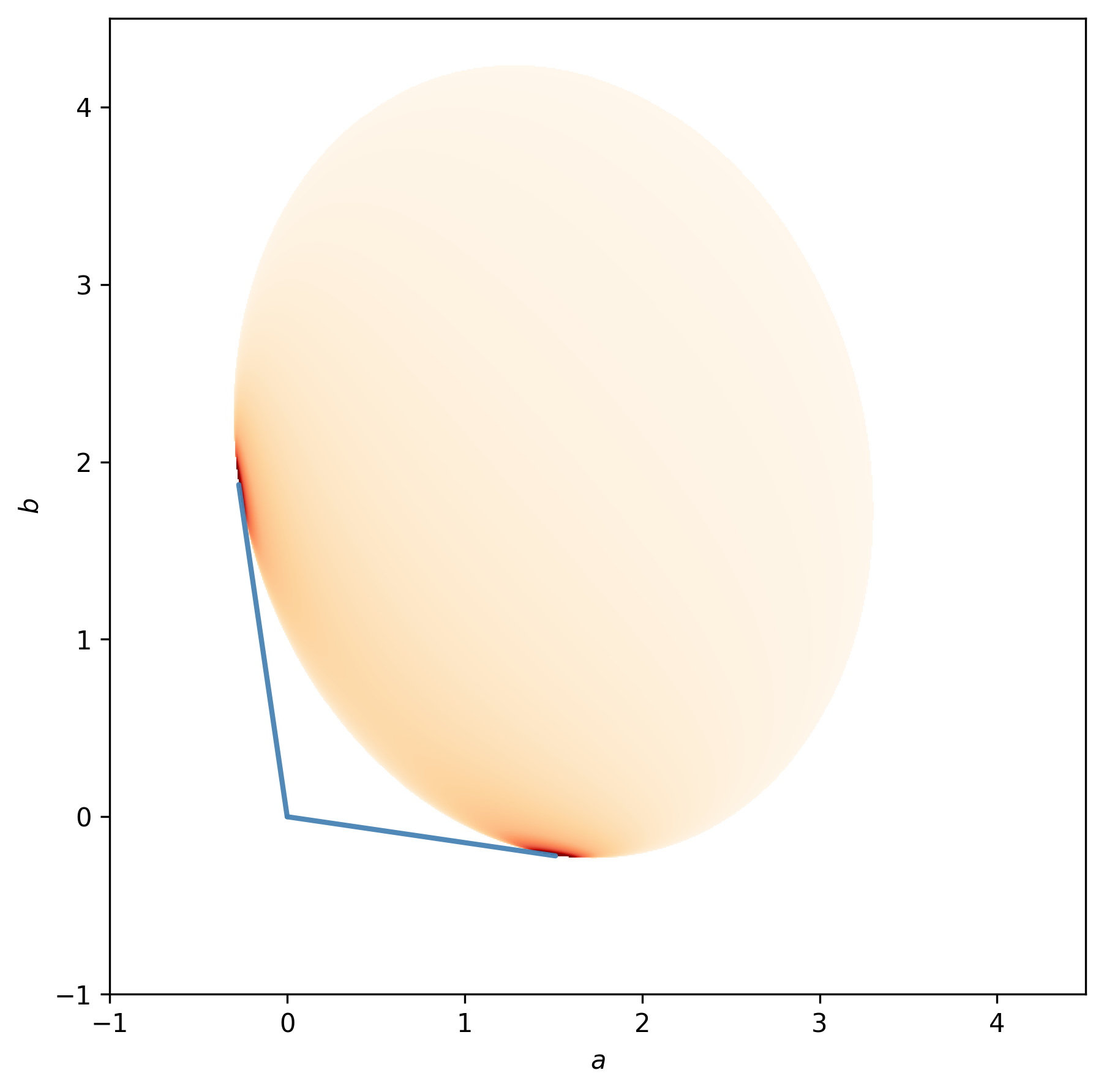}
        \caption{
        $A = \begin{bmatrix} 
            3 & -1 \\ 
            -1 & 0
        \end{bmatrix},
        B = \begin{bmatrix} 
            3 & 2 \\ 
            2 & 1 
        \end{bmatrix}$}
    \end{subfigure}
    \hfill
    \begin{subfigure}[t]{0.48\textwidth}
        \centering
        \includegraphics[width=\linewidth]{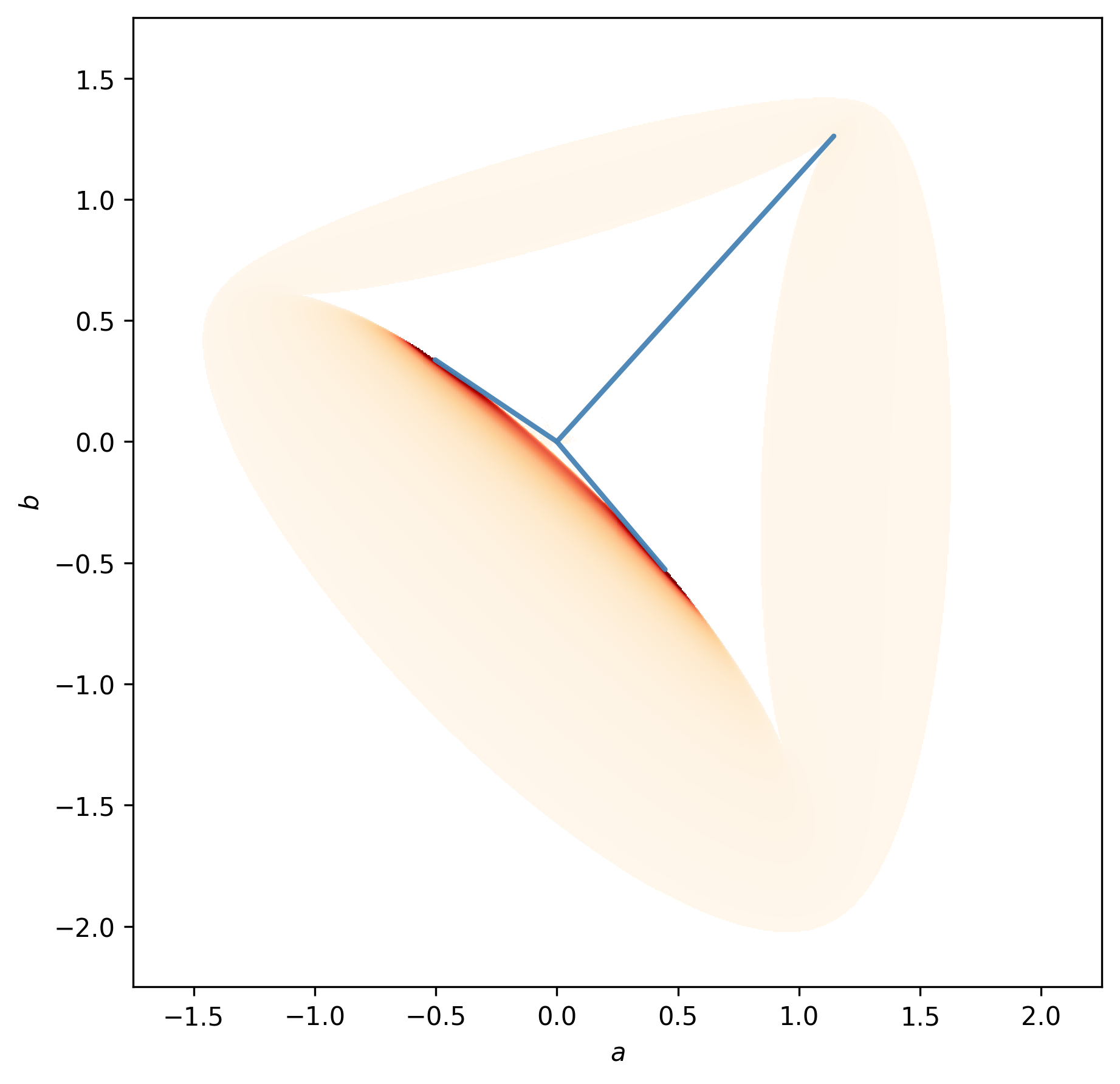}
        \caption{
        $A = \begin{bmatrix} 
            1 & 0.5 & -1 \\ 
            0.5 & 0 & 1 \\
            -1 & 1 & 0 
        \end{bmatrix},
        B = \begin{bmatrix} 
            1 & 1 & 0.5 \\ 
            1 & -1 & -1 \\
            0.5 & -1 & 0
        \end{bmatrix}$}
    \end{subfigure}
    \caption{Sample figures of $\mu_c$ and $\mu_s$. The orange part is the support of $\mu_c$ with color intensity representing density, and the blue segments are the location of the singular part $\mu_s$.}
    \label{fig:sample_figures}
\end{figure}

If $A$ and $B$ commute, the matrix $(I-A/a)(bA-aB)^{-1}$ is Hermitian and consequently the density is $0$ everywhere. In particular, this happens for real diagonal $A$ and $B$. We perturb this trivial diagonal case. Consider
\begin{align*}
    A(\varepsilon) = \underbrace{A(0)}_{\tx{diag}} + \varepsilon A'(0), \, B(\varepsilon) = \underbrace{B(0)}_{\tx{diag}} + \varepsilon B'(0)
\end{align*}
where $A(0)$ and $B(0)$ are real diagonal; $A'(0)$ and $B'(0)$ are general Hermitian matrices. Denote the density of $\mu_c$ by $\rho$ under this setup:
\begin{align}\label{eq:density_in_intro}
    \rho(a,b,\varepsilon) := \frac{1}{2\pi} \sum_{i=1}^n \left|\Im\left( \lambda_i\left(\left(I-\frac{A(\varepsilon)}{a}\right)(bA(\varepsilon)-aB(\varepsilon))^{-1}\right)\right)\right|.
\end{align}

For small $\varepsilon >0$, we expect the density of $\mu_c$ to be small and mostly zero, as illustrated in Figure~\ref{fig:sample_figures_nearly_diagonal} below.

\begin{example}
We consider the following nearly diagonal matrices:
\begin{equation*}
\begin{split}
    &A(\varepsilon) = \begin{bmatrix}
    -0.5 & 0 & 0 & 0 \\
    0 & 0.5 & 0 & 0 \\
    0 & 0 & 0.5 & 0 \\
    0 & 0 & 0 & 1.5
    \end{bmatrix} \;\,\,\,+ \varepsilon \begin{bmatrix}
    -1 & -1 & -1 & 1 \\
    -1 & 1 & 1 & -2 \\
    -1 & 1 & 0 & 1 \\
    1 & -2 & 1 & 1
    \end{bmatrix}, \\
    &B(\varepsilon) = \begin{bmatrix}
    -0.5 & 0 & 0 & 0 \\
    0 & 1.5 & 0 & 0 \\
    0 & 0 & 0.5 & 0 \\
    0 & 0 & 0 & -0.5
    \end{bmatrix} + \varepsilon \begin{bmatrix}
    -1 & 1 & -1 & 1 \\
    1 & -1 & 1 & -2 \\
    -1 & 1 & 0 & 2 \\
    1 & -2 & 2 & 1
    \end{bmatrix}.
\end{split}
\end{equation*}
\label{ex:sample_44}
\end{example}

\begin{figure}[H]
    \begin{subfigure}[H]{0.48\textwidth}
        \centering
        \includegraphics[width=\linewidth]{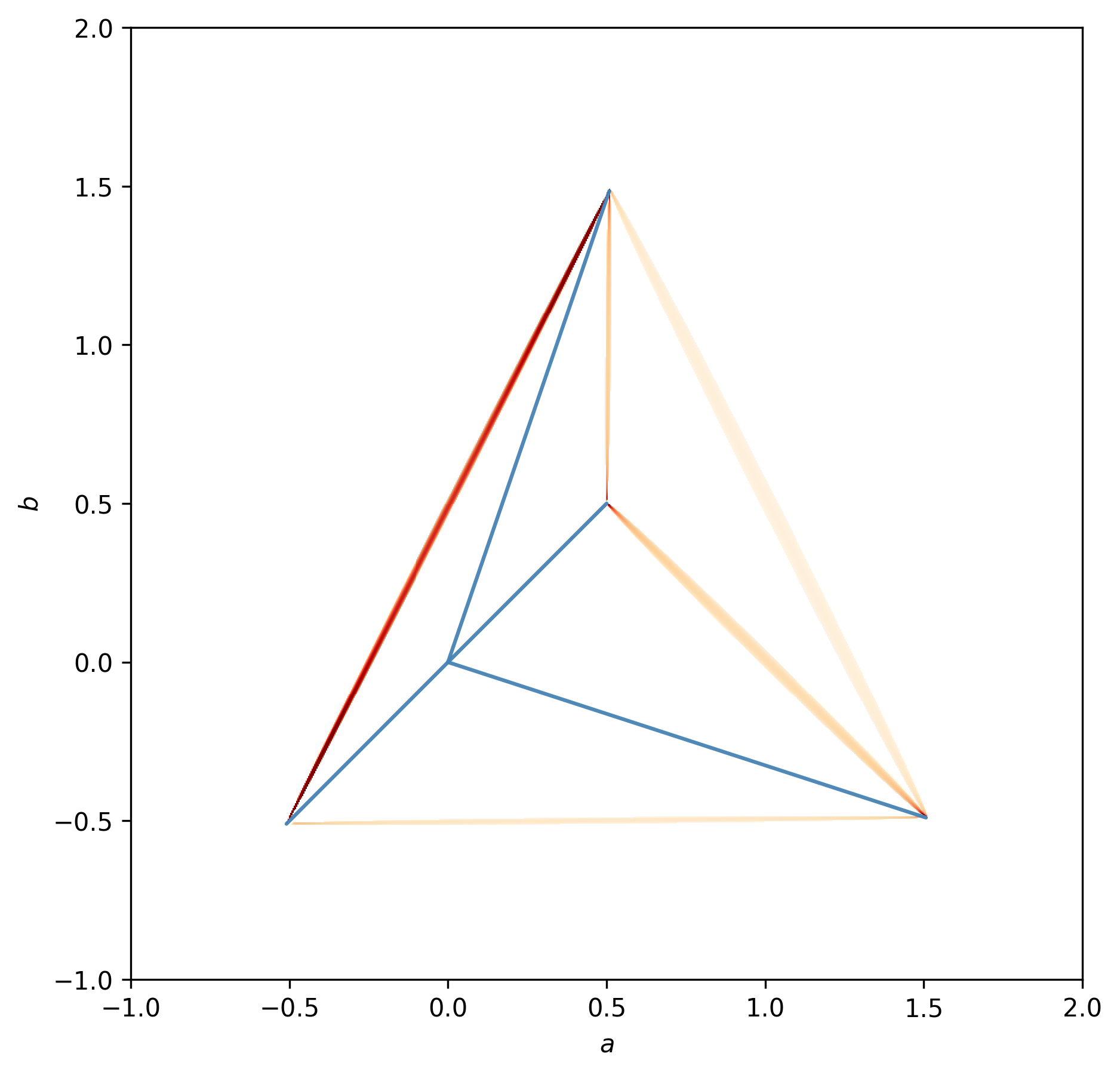}
        \caption{$A = A(0) + 0.01 A'(0)$, $B = B(0) + 0.01B'(0)$}
    \end{subfigure}
    \hfill
    \begin{subfigure}[H]{0.48\textwidth}
        \centering
        \includegraphics[width=\linewidth]{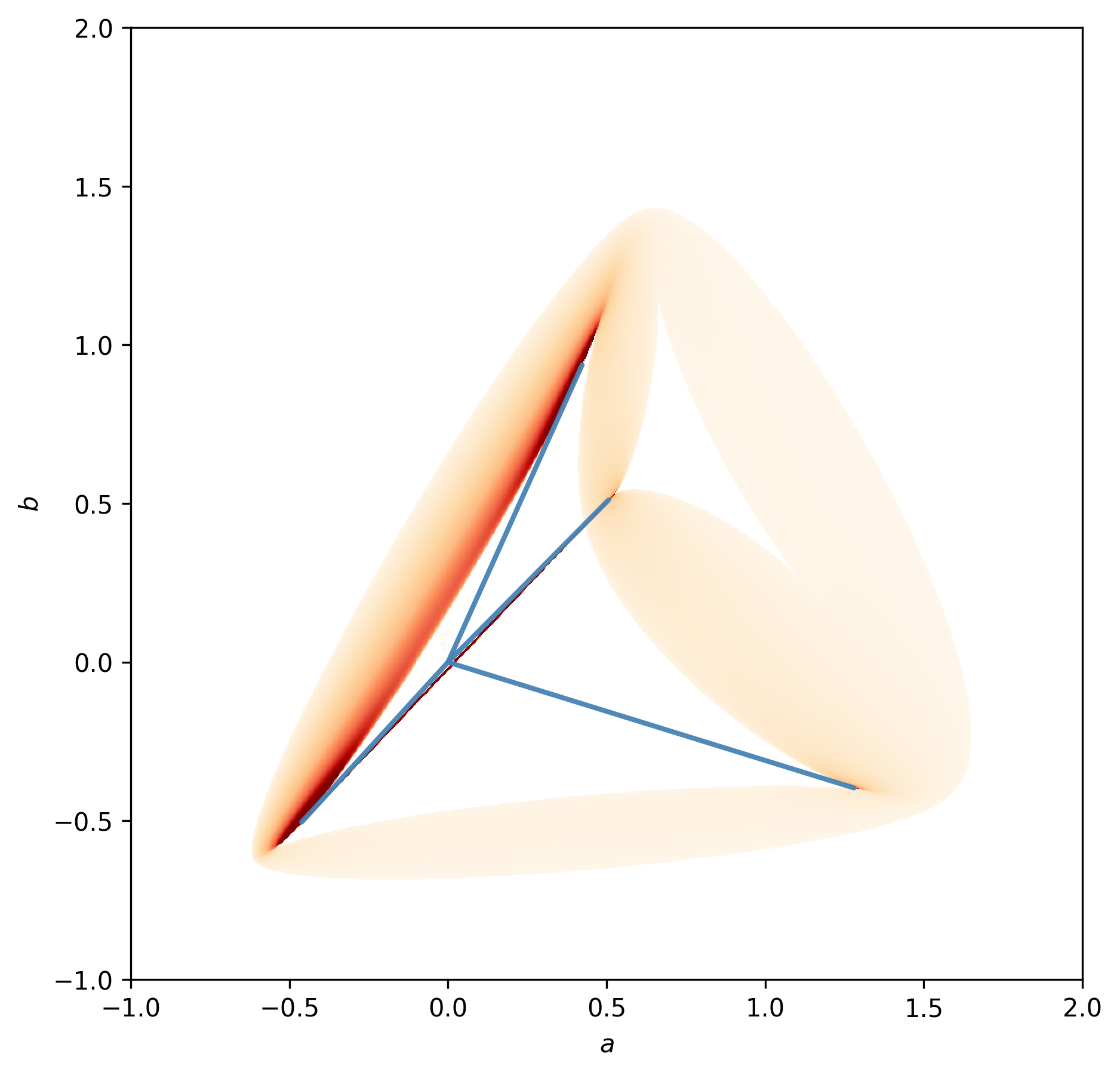}
        \caption{$A = A(0) + 0.1 A'(0)$, $B = B(0) + 0.1B'(0)$}
    \end{subfigure}
    \caption{$\mu_{A,B}$ in a nearly diagonal setup, where the choices of $A(0)$, $B(0)$, $A'(0)$, $B'(0)$ are given in Example~\ref{ex:sample_44}. Again, the orange part is the support of $\mu_c$ and the blue segments are the location of $\mu_s$.}
    \label{fig:sample_figures_nearly_diagonal}
\end{figure}

\begin{figure}[H]
    \centering
    \includegraphics[width=0.55\linewidth]{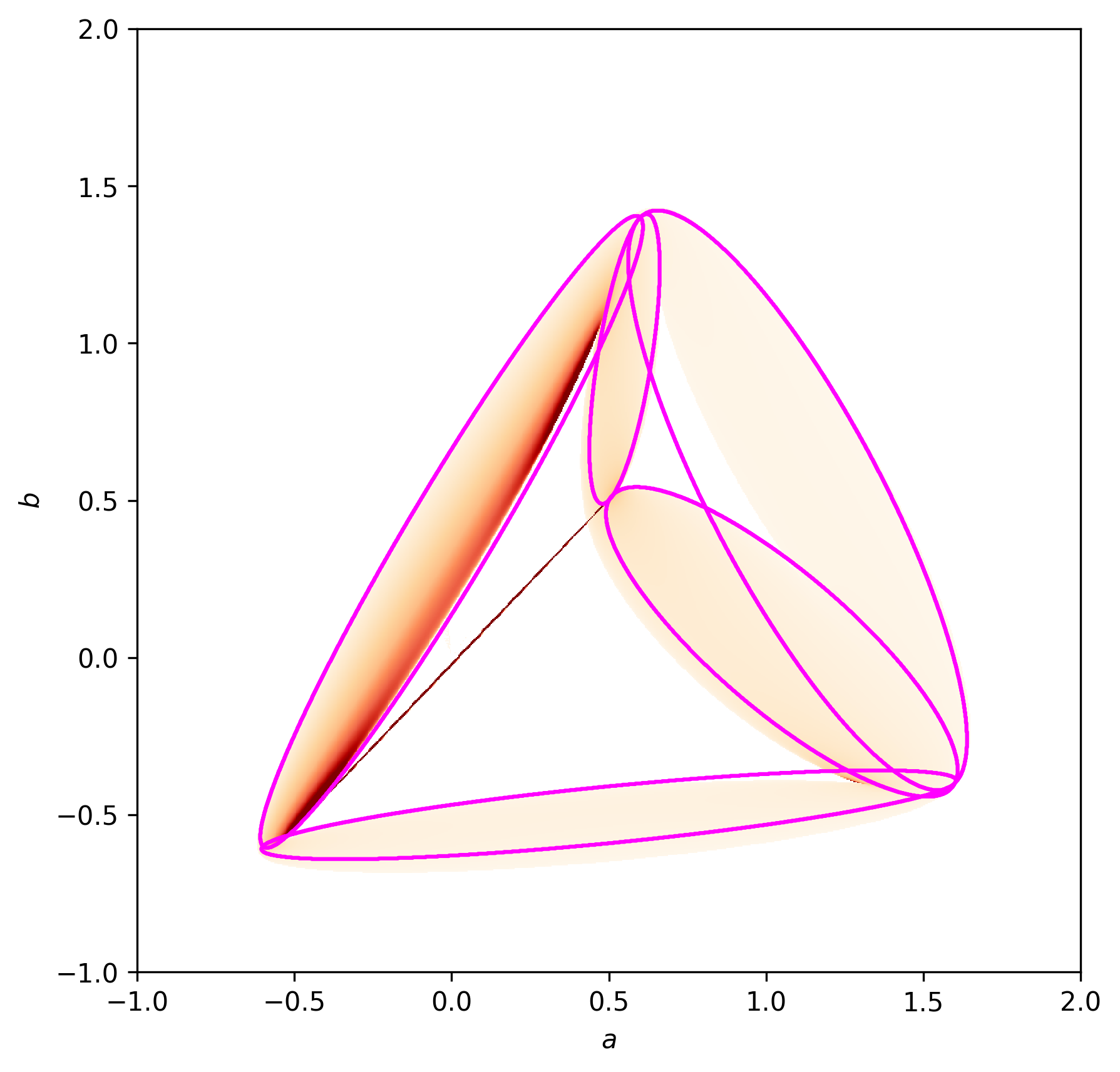}
    \caption{The boundary of $\tx{supp}(\mu_c)$ can be approximated be a union of ellipses.}
    \label{fig:44 with ellipses}
\end{figure}

We focus on the absolutely continuous part $\mu_c$. It appears that the boundary of $\tx{supp}(\mu_c)$ can be approximated by a union of ellipses. We plot one choice of the ellipses in Figure~\ref{fig:44 with ellipses}.

The ellipses are computed using what we call a ``$2\times2$ approximating density'', $\rho_{approx}$. Our paper is centered around this approximating density. The formula of $\rho_{approx}$ is given by

\begin{equation}
\label{eq:density_approx_in_intro}
\begin{split}
    &\rho_{approx}(a,b,\varepsilon) := \sum_{1<i<j<n} \rho_{i,j}(a,b,\varepsilon),\\
    &\rho_{i,j}(a,b,\varepsilon) := \frac{1}{2\pi} \sum_{k=1}^2\left|\Im\,\lambda_k\left(\left(I - \frac{A_{\{i,j\}}(\varepsilon)}{a}\right)\left(bA_{\{i,j\}}(\varepsilon)-aB_{\{i,j\}}(\varepsilon)\right)^{-1}\right)\right|,
\end{split}
\end{equation}

where $A_{\{i,j\}}$ and $B_{\{i,j\}}$ are $2\times2$ principal submatrices of $A$ and $B$ corresponding to the rows/columns $i$ and $j$. We show that $\rho_{approx}$ approximates the true density up to $3/2$ order:

\begin{theorem*}[$2\times2$ Approximation of $\mu_c$, Simplified Version of Theorem~\ref{thm:22_approx}] There exists some $\varepsilon^*>0$ and $C>0$ such that for every $(a,b)\in \mathbb R^2$ outside a small set and every $\varepsilon\in [0,\varepsilon^*]$, we have
\begin{align*}
    \left|\rho(a,b,\varepsilon) - \rho_{approx}(a,b,\varepsilon)\right| \leq C \varepsilon^{3/2},
\end{align*}
where $\rho$ and $\rho_{approx}$ are as in \eqref{eq:density_in_intro} and \eqref{eq:density_approx_in_intro}.
\end{theorem*}

\begin{figure}[H]
    \centering
    \includegraphics[width=0.55\linewidth]{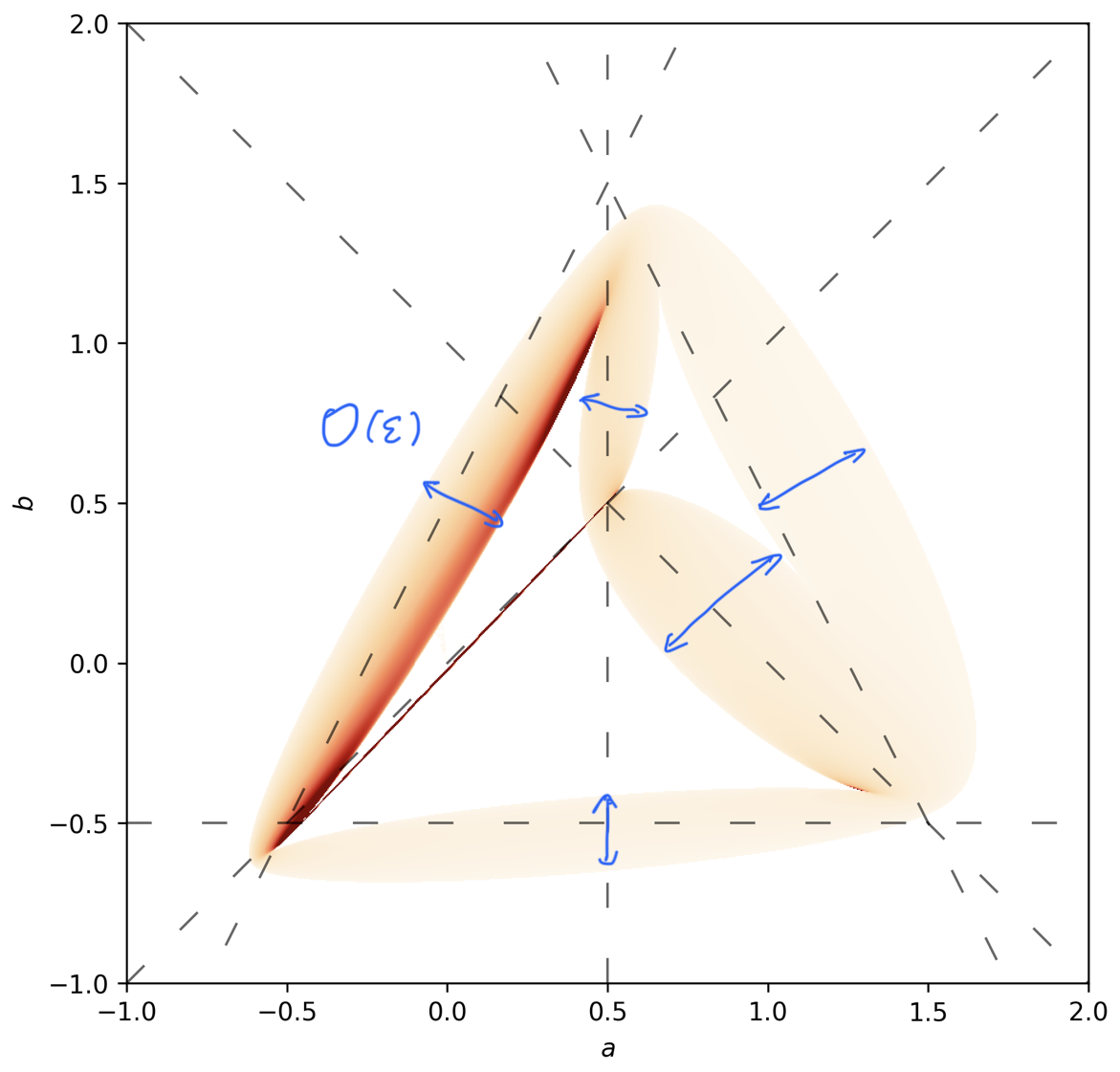}
    \caption{The support lies around collision lines (dashed lines), and locally around each collision line, the perpendicular width (blue arrowed segments) of the support is $O(\varepsilon)$.}
    \label{fig:width_support_intro}
\end{figure}

Additionally, we use $\rho_{approx}$ to study the geometric properties of the true density. We prove that the support of the $\mu_c$ lie around certain ``collision lines'' on the plane; and locally near each collision line, its width measured perpendicular to the line is $O(\varepsilon)$, as illustrated in Figure~\ref{fig:width_support_intro}. The density of $\mu_c$ is $O(\varepsilon)$ as well. 

\begin{remark}
For every pair of commuting Hermitian matrices $A$ and $B$, there exist real diagonal $A_0$ and $B_0$ such that the $\mu_c$ of the two pairs are equal. By~\cite[Theorem~2.5.5]{HornJohnson2013}, there exists a unitary matrix $U$ such that $A_0 := U^*A U$ and $B_0 := U^*BU$ are real diagonal, and $\mu_{A,B} = \mu_{A_0, B_0}$ by the unitary invariance of the tracial joint spectral measures~\cite[Proposition~3.8]{Hei_spectral}. Hence, to study the behavior of $\mu_c$ for general nearly commuting matrices, it suffices to consider nearly diagonal ones.
\end{remark}

\section{Setup and Preliminaries}

\subsection{Setup and Assumptions}
\label{sec:setup_assumptions}

Consider affine Hermitian matrix-valued functions $A(\varepsilon) = A(0) + \varepsilon A'(0)$ and $B(\varepsilon) = B(0) + \varepsilon B'(0)$, where $A(0)$ and $B(0)$ are furthermore real diagonal. Define matrix $M$ by
\begin{align*}
    &M : \mathbb R^3\, \backslash \{(a,b,\varepsilon) \tx{ where } a = 0 \tx{ or } bA(\varepsilon)-aB(\varepsilon) \tx{ not invertible}\} \to M_n(\mathbb C) \\
    &M(a,b,\varepsilon):=\left(I-\frac{A(\varepsilon)}{a}\right)(bA(\varepsilon)-aB(\varepsilon))^{-1}.
\end{align*}

We discuss the density of the absolutely continuous part of the tracial joint spectral measure (Theorem~\ref{thm:density_of_measure_original}):
\begin{align*}
    \rho(a,b,\varepsilon) &:= \frac{1}{2\pi} \sum_{i=1}^n \left|\Im\left(\lambda_i\left(\left(I - \frac{A(\varepsilon)}{a}\right)(bA(\varepsilon)-aB(\varepsilon))^{-1}\right)\right)\right|\\
    &= \frac{1}{2\pi} \sum_{i=1}^n|\Im(\lambda_i(M(a,b,\varepsilon)))|.
\end{align*}
At $\varepsilon = 0$, $M$ is real diagonal $\rho(\cdot, \cdot, 0) = 0$ everywhere. We study the eigenvalues for small $\varepsilon>0$.\\

We make the following assumption throughout the paper.

\begin{assumption}
\label{asmp:non_collinear}
Consider points $(A(0)_{i,i}, B(0)_{i,i})$ determined by diagonal entries of $A(0)$ and $B(0)$. We assume no three such points are collinear. Equivalently, we assume for any $i<j<k\in[1,n]$,
\begin{align*}
    (A(0)_{i,i} - A(0)_{j,j})(B(0)_{i,i}-B(0)_{k,k}) - (A(0)_{i,i} - A(0)_{k,k})(B(0)_{i,i}-B(0)_{j,j}) \neq 0.
\end{align*}
\end{assumption}

This assumption ensures certain ``collision lines'', defined in Definition~\ref{def:collision_lines} and shown as dashed lines in Figure~\ref{fig:width_support_intro}, are distinct.

\subsection{Preliminary Properties of $\rho$ and $M$}
\label{sec:preliminary_properties_rho_M}

\begin{proposition}[{\cite[Proposition~4]{Hei_spectral_thesis}}]
For any $\varepsilon_0\geq 0$, $\tx{supp}(\rho(\cdot,\cdot,\varepsilon_0))$ is compact.
\end{proposition}

\begin{definition}
\label{def:Sigma}
Denote the set where $M$ is not well-defined by $\Sigma$ and further overload the notation to refer to the slices of $\Sigma$:
\begin{align*}
    \Sigma &:= \{(a,b,\varepsilon)\,|\, a = 0 \,\tx{or}\, \det(bA(\varepsilon) - aB(\varepsilon)) = 0\},\\
    \Sigma(\varepsilon_0) &:= \{(a,b) \,|\, (a,b,\varepsilon_0) \in \Sigma\},\\
    \Sigma(a_0,b_0) &:= \{\varepsilon \,|\, (a_0, b_0, \varepsilon) \in \Sigma\}.
\end{align*}
\end{definition}

\begin{proposition}[Global Well-Definedness of $M$]
\label{prop:Sigma_finite_lines}
\leavevmode
\begin{enumerate}
    \item For every $\varepsilon$, $\Sigma(\varepsilon)$ is either the entire $\mathbb R^2$ or a finite union of lines.
    \item Moreover, if $\Sigma \neq \mathbb R^3$, $\Sigma(\varepsilon)$ is a finite union of lines for all $\varepsilon$ but finitely many.
    \item For every $(a, b)$, $\Sigma(a,b)$ is either the entire $\mathbb R$ or finitely many points.
    \item Moreover, if $\Sigma\neq \mathbb R^3$, $\Sigma(a, b)$ is finitely many points for all $(a,b)$ but those belong to a finite union of lines.
\end{enumerate}
\end{proposition}
\begin{proof}
First we prove (1) and (3). Assuming $a \neq 0$, we notice:
\begin{align*}
    \det(bA(\varepsilon) - aB(\varepsilon)) = 0 &\iff a^n \det(\frac{b}{a}A(\varepsilon) - B(\varepsilon)) = 0\\
    &\iff \det(\frac{b}{a}A(\varepsilon) - B(\varepsilon)) = 0.
\end{align*}
Fixing $\varepsilon = \varepsilon_0$, the last determinant becomes a polynomial in $b/a$; and it is either the zero polynomial, in which case $\Sigma(\varepsilon_0) = \mathbb R^2$, or has finitely many roots, in which case $\Sigma(\varepsilon_0)$ is a finite union of lines. This proves statement (1). \\

Alternatively, fixing $a=a_0,\,b=b_0$ turns the determinant into a polynomial in $\varepsilon$, which is either zero for every $\varepsilon \in \mathbb R$ or has finitely many roots. This proves statement (3).\\

Now we prove statement (2). Assume $\Sigma\neq \mathbb R^3$ and pick some $(a_0, b_0, \varepsilon_0)\notin \Sigma$. We know $\det(b_0A(\varepsilon)- a_0B(\varepsilon))$ is a polynomial in $\varepsilon$ and is not the zero polynomial, since it is non-zero at $\varepsilon_0$. Consequently $\det(b_0A(\varepsilon)- a_0B(\varepsilon))$ has only finitely many roots. Pick some $\varepsilon_1$ that is not a root, we show $\Sigma(\varepsilon_1)$ is a finite union of lines. \\

We know $\Sigma(\varepsilon_1) = \{a=0\} \cup \{(a,ta) \tx{ where } t \tx{ s.t. } \det(tA(\varepsilon_1) - B(\varepsilon_1)=0\}$. We show the second component is finitely many lines, equivalently that $\det(tA(\varepsilon_1) - B(\varepsilon_1))$ is not the zero polynomial. And it is indeed not the zero polynomial because $t = b_0/a_0$ is not a root.\\

Similar argument proves statement (4).
\end{proof}

\begin{proposition}[Local Well-Definedness of $M$]
\label{proposition:M_well_defined_locally}
Assume $M$ is well-defined at some $(a_0, b_0, \varepsilon_0)$. There exists a neighborhood $U\subseteq \mathbb R^3$ around $(a_0, b_0, \varepsilon_0)$ such that $M$ is well-defined at any point in $U$. 
\end{proposition}
\begin{proof}
$\{(a,b,\varepsilon) \,|\, a \neq 0 \tx{ and } \det(bA(\varepsilon) - a B(\varepsilon)) \neq 0\}$ is an open set.
\end{proof}

\begin{proposition}[Eigenvalues Form Conjugate Pairs]\label{prop:eval_form_conj_pairs}
For any $(a,b,\varepsilon)$ where $M$ is well-defined, its eigenvalues are either real or form conjugate pairs.
\end{proposition}
\begin{proof}
The eigenvalues of $M$ are the roots of the real polynomial
\begin{align*}
    \det(\lambda I - M) = \det(bA-aB)^{-1} \det \left(\lambda (bA-aB) - \left(I- \frac{A}{a}\right) \right).
\end{align*}
\end{proof}

\section{Local Separation of Eigenvalues}
\label{sec:local_sep_evals}

The goal of this section is to find a suitable environment where the density $\rho$ could be computed more easily. We shall conclude with the following theorem in this section:

\begin{theorem}
[Theorem~\ref{thm:separation_eval_near_collision_line}, Simplified Form] For any $L_{i,j}$ (as in Definition~\ref{def:collision_lines}) and almost every $(a_0, b_0)\in L_{i,j}$, there exists a neighborhood $U$ around $(a_0, b_0, 0)$ and $2\times2$ matrix $T$ defined on $U$ such that 
\begin{align*}
    \rho(a,b,\varepsilon) = \frac{1}{2\pi} (|\Im\,\lambda_1(T(a,b,\varepsilon))| + |\Im\,\lambda_2(T(a,b,\varepsilon))|)
\end{align*}
everywhere in $U$. Furthermore, $T$ is approximated entry-wise in $U$ by
\begin{align*}
    M^{\{i,j\}}(a,b,\varepsilon) := \left(I - \frac{A_{\{i,j\}}}{a}\right)(bA_{\{i,j\}} - aB_{\{i,j\}})^{-1}.
\end{align*}
\end{theorem}

We pack up the objects and statements into one compact form, the Local Separation of Eigenvalues setup, or the LSE setup, defined in Definition~\ref{def:LSE}.

\subsection{Properties of Eigenvalues Near a Collision Line}
\label{sec:properties_eval_near_line}

We begin with a definition:

\begin{definition}
\label{def:collision_lines}
(Collision Lines) We define collision lines $L_{i,j}$ ($1\leq i<j\leq n$) as
\begin{align*}
    L_{i,j} :=\{(A(0)_{i,i}+t(A(0)_{j,j}-A(0)_{i,i}), B(0)_{i,i}+t(B(0)_{j,j}-B(0)_{i,i})),\, t \in \mathbb{R}\} \backslash \Sigma
\end{align*}
where $\Sigma$ is defined in Definition~\ref{def:Sigma} and denotes the set of $(a,b,\varepsilon)$ where $M$ is not well-defined. 
\end{definition}
Figure~\ref{fig:44_collision_lines} illustrates the collision lines as dashed lines\\

We note that the non-collinear assumption, Assumption~\ref{asmp:non_collinear}, ensures (and is in fact equivalent to) that the collision lines are distinct. We discuss properties of collision lines.

\begin{proposition}
For any $(a,b)\in \mathbb R^2$ where $M(a,b,0)$ is well-defined, two eigenvalues of $M(a,b,0)$ are equal if and only if $(a,b)$ belongs to some collision line. Moreover, 
\begin{align*}
    (M(a,b,0))_{i,i} = (M(a,b,0))_{j,j} \iff (a,b) \in L_{i,j}.
\end{align*}
\end{proposition}
\begin{proof}
Notice that the matrix $M(a,b,0)$ is diagonal:
\begin{align*}
    M(a,b,0) &= \left(I - \frac{A(0)}{a}\right)(bA(0) - aB(0))^{-1}\\
    &= \tx{diag}(\frac{1-A(0)_{1,1}/a}{bA(0)_{1,1}-aB(0)_{1,1}}, ..., \frac{1-A(0)_{n,n}/a}{bA(0)_{n,n}-aB(0)_{n,n}}).
\end{align*}
Denote the eigenvalues of $M(a,b,0)$ by $\lambda_i(a,b,0)$, in the order that they appear in the diagonal entries of $M$:
\begin{align*}
    \lambda_i(a,b,0) &:= (M(a,b,0))_{i,i} = \frac{1-A(0)_{i,i}/a}{bA(0)_{i,i}-aB(0)_{i,i}} \in \mathbb R.
\end{align*}
Set two eigenvalues equal and solve for $(a,b)$. With a little bit computation, we observe:
\begin{align*}
    &\lambda_i(a,b,0) = \lambda_j(a,b,0) \\
    \iff& \frac{1-A(0)_{i,i}/a}{bA(0)_{i,i}-aB(0)_{i,i}} = \frac{1-A(0)_{j,j}/a}{bA(0)_{j,j}-aB(0)_{j,j}}\\
    \iff& b(A(0)_{j,j}-A(0)_{i,i}) -a(B(0)_{j,j}-B(0)_{i,i}) +A(0)_{i,i}B(0)_{j,j} -A(0)_{j,j}B(0)_{i,i}=0\\
    \iff& (a,b) \in L_{i,j}.
\end{align*}
\end{proof}

More than two eigenvalues become equal on the intersection of the collision lines. For example, we have
\begin{align*}
    \lambda_1(a,b,0) = \lambda_2(a,b,0) = \lambda_3(a,b,0) \iff (a,b) \in L_{1,2}\cap L_{1,3} \cap L_{2,3}.
\end{align*}

We go back to the $4\times4$ perturbation example, Example~\ref{ex:sample_44}, and plot the collision lines:
\begin{figure}[H]
    \centering
    \includegraphics[width=0.476\linewidth]{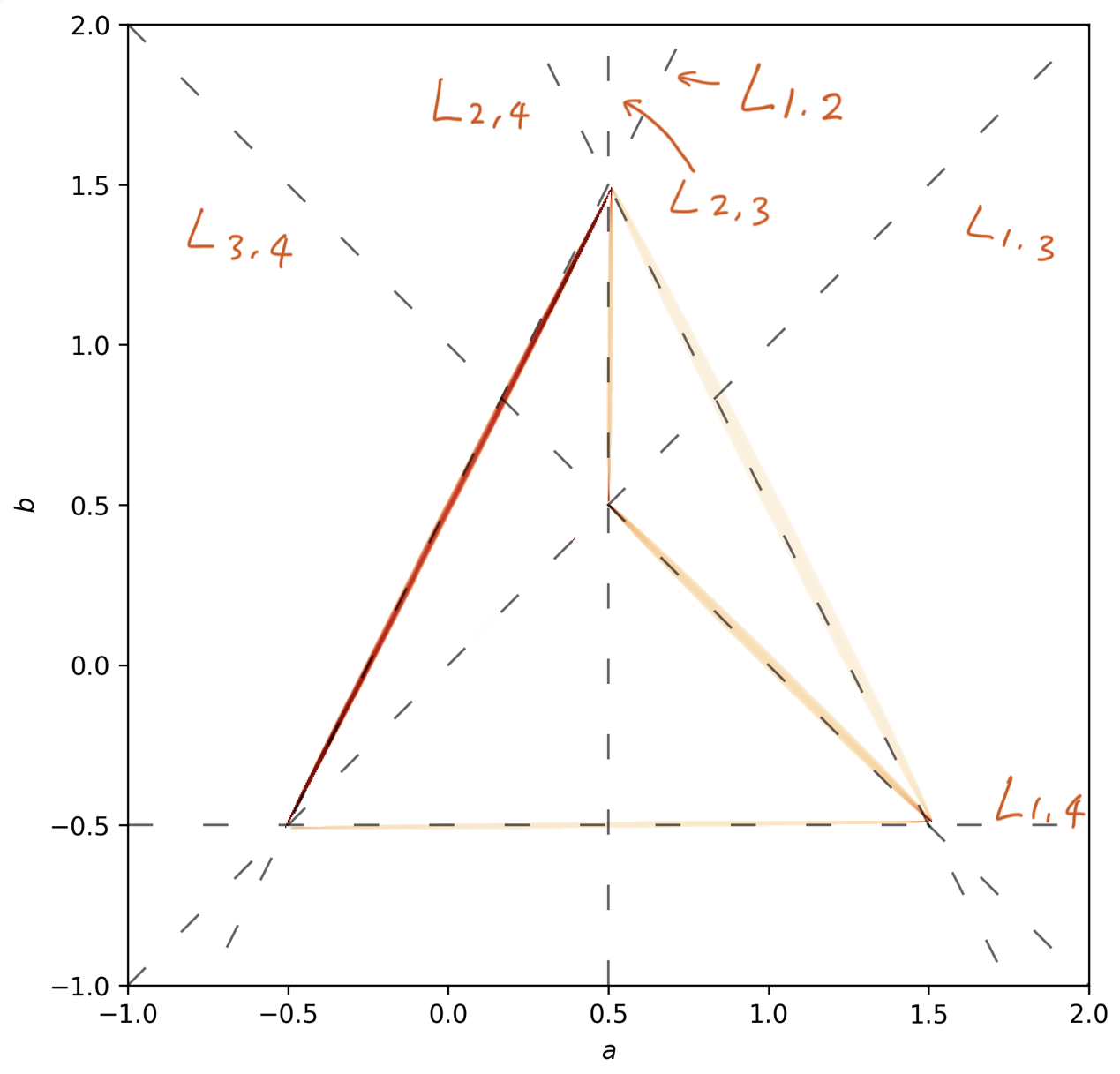}
    \hfill
    \includegraphics[width=0.484\linewidth]{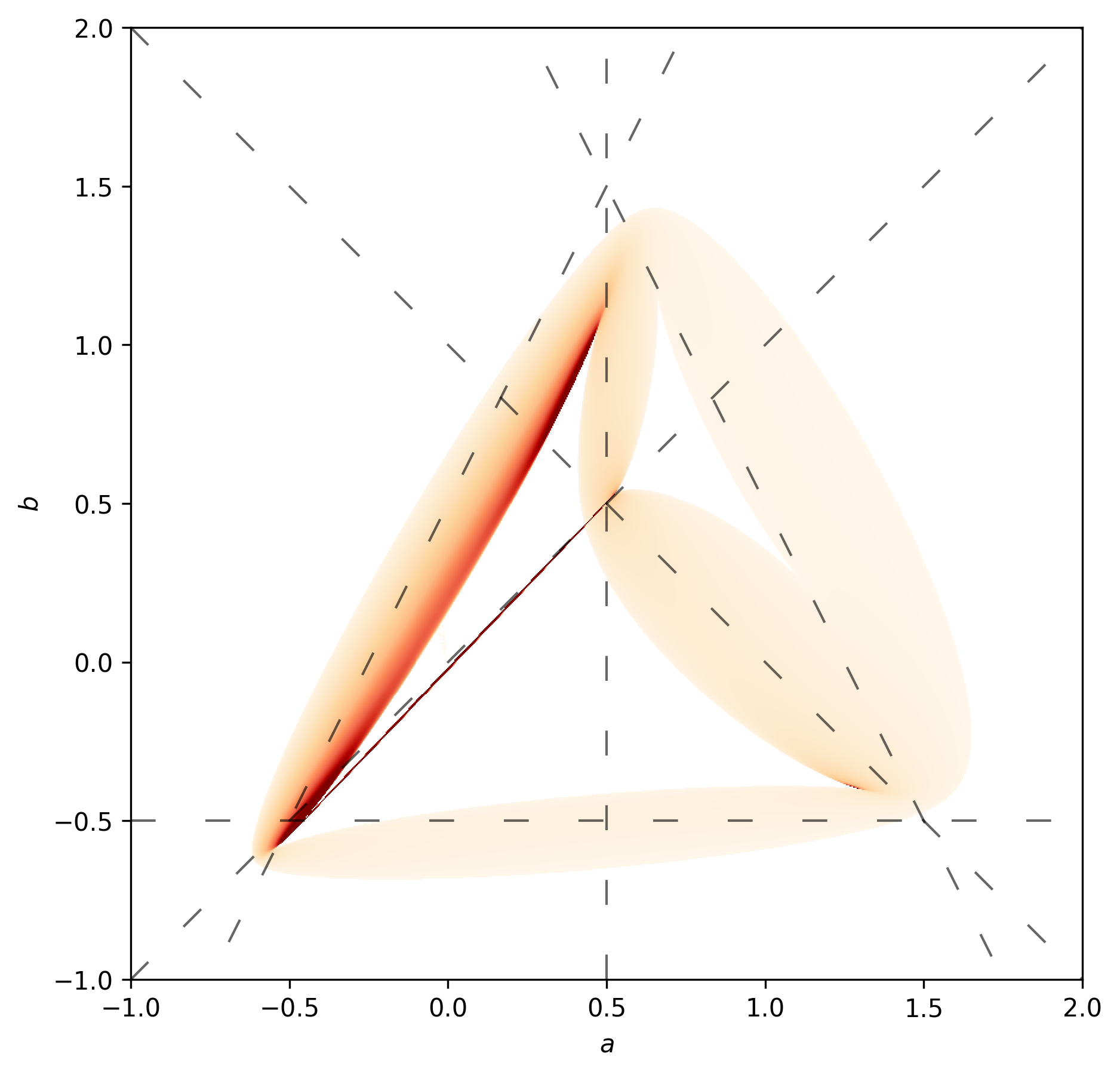}
    \caption{Collision lines in the case of Example~\ref{ex:sample_44}. The background density is plotted using $\varepsilon = 0.01$ (left) and $\varepsilon = 0.1$ (right).}
    \label{fig:44_collision_lines}
\end{figure}

We observe that for small $\varepsilon$, the density is positive at some $(a,b)$ only if $(a,b)$ sits near a collision line. 

\begin{theorem}
\label{thm:rho_locally_0_off_collision_line}
Consider $(a_0, b_0, 0)$ that does not belong to any collision line. There exists a neighborhood $U$ around $(a_0, b_0, 0)$ such that $\rho(a,b,\varepsilon)=0$ everywhere in $U$.
\end{theorem}
\begin{proof}
Put an open ball around each real simple eigenvalue of $M(a_0,b_0,0)$ and denote the open balls by $\{B_i\}$, $i = 1,2,...,n$. By Lemma~\ref{lemma:kato_local_separation_of_eigenvalues}, the eigenvalues of a matrix depend continuously on the matrix when viewed as an unordered set. So there exists a neighborhood $U$ around $(a_0,b_0,0)$ such that for every $(a,b,\varepsilon) \in U$, $M(a,b,\varepsilon)$ has precisely one eigenvalue inside each $B_i$, and no conjugate pairs can be formed. On the other hand, we know by Proposition~\ref{prop:eval_form_conj_pairs} that the eigenvalues of $M(a,b,\varepsilon)$ is either real or form conjugate pairs, so they must be real. Hence $\rho(a,b,\varepsilon) = 0$.
\end{proof}

\begin{lemma}[{\cite[Theorem II.5.14]{Kato1995}}]\label{lemma:kato_local_separation_of_eigenvalues}
For any fixed matrix $T$, $\tx{dist}(\Lambda(T),\Lambda(T+A)) \to 0$ as $\|A\|\to 0$.
\end{lemma}

Here $\Lambda(T) = \{\lambda_1(T), ..., \lambda_n(T)\}$ denotes the unordered set of eigenvalues of an $n\times n$ matrix $T$, and the distance between two set of eigenvalues $\Lambda(T)$, $\Lambda(S)$ is defined as
\begin{align*}
    \tx{dist}(\Lambda(T), \Lambda(S)) := \min_{\sigma \in S_n}\max_{1\leq i\leq n} |\lambda_i(T) - \lambda_{\sigma(i)}(S)|,
\end{align*}
where $\sigma$ takes value in the permutation group $S_n$.\\

We know more about the simple eigenvalues of $M$. Consider an $(a_0, b_0)$ that belongs to precisely one collision line, say $L_{1,2}$; we know $M(a_0, b_0,0)$ has two repeated eigenvalues and $n-2$ simple ones. By Lemma~\ref{lemma:kato_local_separation_of_eigenvalues}, the $n-2$ eigenvalues should remain simple around $(a_0, b_0, 0)$. By~\cite[p117, footnote 1]{Kato1995}, these simple eigenvalues are in fact analytic.

\begin{lemma}[{\cite[p117, footnote 1]{Kato1995}}]\label{lemma:kato_simple_evals_are_analytic}
The simple eigenvalues of an analytic matrix-valued function are analytic.
\end{lemma}

We formulate a setup and summarize our observations inside this setup.\\
\newline\newline

\begin{definition}[The LSE setup]
\label{def:LSE}
\leavevmode
\begin{itemize}
    \item Consider some collision line, say $L_{1,2}$ and $L_{1,2} = \{(1,b)| b \in \mathbb R\}$; pick a point that belongs to $L_{1,2}$ but not to any other collision lines, say $(1,0)$. Choose $a^*, b^*, \varepsilon^*, r$ such that the statements of Theorem~\ref{thm:local_properties_eval_near_collision_line} hold inside the rectangle $[1-a^*,1+a^*]\times[-b^*,b^*]\times[0,\varepsilon^*]$. We pack up this choice of the collision line, the point and the rectangle, and the conclusions of Theorem~\ref{thm:local_properties_eval_near_collision_line} as the background environment of the Local Separation of Eigenvalues (LSE) setup.
    \item (Theorem~\ref{thm:separation_eval_near_collision_line}) Additionally, there exists a $2\times 2$ matrix $T$ whose eigenvalues equal to the two potentially non-real eigenvalues of $M$, and another $2\times2$ matrix $M^{2\times2}$ that approximates $T$ entry-wise up to second order.
    \item (Theorem~\ref{thm:width_of_support_and_density}) Finally, there exists $\varepsilon_{cvx}^* \in [0,\varepsilon^*]$ such that $\tx{Disc}(\lambda I - T) := \tx{Disc}_{\lambda} \det(\lambda I - T))$ is convex in $[1-a^*,1+a^*]\times[-b^*,b^*]\times[0,\varepsilon_{cvx}^*]$.
\end{itemize}
We prove the three parts gradually in this paper. For easier notation, we refer to the first part, i.e., the background environment, as ``LSE(1)''; the first and the second parts, i.e., the background environment together with the statement of Theorem~\ref{thm:separation_eval_near_collision_line}, as ``LSE(2)''; and refer to the three parts combined as \textbf{the LSE setup}.
\end{definition}

\begin{theorem}[Local Properties of Eigenvalues Near a Collision Line]
\label{thm:local_properties_eval_near_collision_line} 
Consider some collision line, say $L_{1,2}$ and $L_{1,2} = \{(1,b)| b \in \mathbb R\}$; pick a point that belongs to $L_{1,2}$ but not to any other collision lines, say $(1,0)$. There exists positive constants $a^*, b^*, \varepsilon^*, r$ such that the following holds inside the rectangle $[1-a^*,1+a^*]\times[-b^*,b^*]\times[0,\varepsilon^*]$:
\begin{enumerate}
    \item $M(a,b,\varepsilon)$ is well-defined and analytic.
    \item Two eigenvalues of $M(a,b,\varepsilon)$ lie in $B_r(M(1,0,0)_{1,1})$, while the other $n-2$ are strictly outside and simple.
    \item Sort the eigenvalues that outside $B_r(M(1,0,0)_{1,1})$ in decreasing order and label them as $\lambda_3(a,b,\varepsilon), ..., \lambda_n(a,b,\varepsilon)$. Each $\lambda_i(a,b,\cdot)$ is analytic on $[0,\varepsilon^*]$.
\end{enumerate}
\end{theorem}

\begin{proof}
For (1), $M$ is well-defined in the rectangle by Proposition~\ref{proposition:M_well_defined_locally}. $M$ is analytic because it is rational. \\

For (2), we know two eigenvalues of $M(1,0,0)$ lie inside $B_r(M(a,b,0)_{1,1})$. Put disjoint balls $B_i$, $i = 3, ..., n$, around the other $n-2$ eigenvalues of $M(1,0,0)$. By Lemma~\ref{lemma:kato_local_separation_of_eigenvalues}, locally around $(1,0,0)$, there should continue to be two eigenvalues inside $B_r(M(a,b,0)_{1,1})$ and one eigenvalue inside each $B_i$.\\

Statement (3) is given by Lemma~\ref{lemma:kato_simple_evals_are_analytic}.
\end{proof} 

\begin{remark}
We show why we could assume the line $L_{1,2}$ is $\{(1,b)\}$. Say we have a general $A(\varepsilon) = A(0) + \varepsilon A'(0)$ and $B(\varepsilon) = B(0) + \varepsilon B'(0)$. Consider the rotation matrix $R$:
\begin{align*}
    R = \frac{1}{\sqrt{(\Delta a)^2 + (\Delta b)^2}} \begin{bmatrix}
        \Delta b & -\Delta a\\
        \Delta a & \Delta b
    \end{bmatrix},
\end{align*}
where $\Delta a := A_{1,1} - A_{2,2}$, $\Delta b := B_{1,1} - B_{2,2}$. We apply $R$ to $A$ and $B$: $C(\varepsilon)=R_{11}A(\varepsilon) + R_{12}B(\varepsilon)$, $D(\varepsilon) = R_{21}A(\varepsilon) + R_{22}B(\varepsilon)$. This rotates the collision line to a vertical position, say $\{(c_0, b)\}$, for some constant $c_0$. Then we consider the combined transformation $V(A,B) = ((R_{11}A + R_{12} B)/c_0,\,c_0 (R_{21}A + R_{22} B))$. This transformation is linear and invertible; so by Proposition~3.8 of Hein{\"a}vaara's paper~\cite{Hei_spectral}, we have $\mu_{\tilde C,\tilde D}$ given by the push forward measure $V_*(\mu_{A,B})$, and hence the densities have a simple relationship.
\end{remark}

\subsection{Separation of Eigenvalues Near a Collision Line}
\label{sec:sep_eval_near_line}

Consider the LSE(1) setup. For each $(a,b,\varepsilon)$ in the rectangle, let $P(a,b,\varepsilon)$ be the projection to the generalized eigenspace of the two eigenvalues of $M(a,b,\varepsilon)$ that belong to the ball $B_r(M(1,0,0))$, which is analytic by \cite[p117, footnote 1]{Kato1995}. Construct basis $S$:
\begin{align*}
    &S(a,b,\varepsilon): [1-a^*, 1+a^*]\times[-b^*, b^*]\times[0,\varepsilon^*] \to M_n(\mathbb C)\\
    &S(a,b,\varepsilon) := [P e_1, P e_2, (I-P) e_{3}, ..., (I-P) e_n],
\end{align*}
where $\{e_1, ..., e_n\}$ is the standard basis for $\mathbb C^n$.\\

We know $P(a,b,0) = I_2 \oplus 0_{n-2}$ for any $(a,b)$. So by continuity, $P(a, b, \varepsilon)$ should stay close to $I_2 \oplus 0_{n-2}$ for sufficiently small $\varepsilon$. Since $S(a, b, 0) = I_n$, it should also stay close to $I_n$ for sufficiently small $\varepsilon$. The following statements formulate this intuition:

\begin{lemma}
\label{lemma:S_error}
Consider the LSE(1) setup. View $S$ in a $(2\times2)\oplus ((n-2)\times(n-2))$ block structure. There exists $C>0$ such that for every $(a,b,\varepsilon)\in[1-a^*, 1+a^*]\times[-b^*, b^*]\times [0,\varepsilon^*]$ and every entry $(i,j)$ of the respective blocks, we have
\begin{align*}
    &|S(a,b,\varepsilon)_{ij} - \mathbf{1}_{i = j}| \leq C\varepsilon^2, \quad |S^{-1}(a,b,\varepsilon)_{ij} - \mathbf{1}_{i = j}| \leq C\varepsilon^2, \tag{for diag blocks}\\
    &|S(a,b,\varepsilon)_{ij}| \leq C\varepsilon, \quad |S^{-1}(a,b,\varepsilon)_{ij}| \leq C\varepsilon. \tag{for off-diag blocks}
\end{align*}
In other words, we may write:
\begin{align*}
    S(u,b,\varepsilon), S^{-1}(u,b,\varepsilon) = \begin{bmatrix}
        I_2 + O(\varepsilon^2) & O(\varepsilon)\\
        O(\varepsilon) & I_{n-2} + O(\varepsilon^2)
    \end{bmatrix}.
\end{align*}
\end{lemma}
\begin{proof}
First we claim $P$ satisfy a similar error bound, that is, for any $(a,b,\varepsilon)$, we have:
\begin{align}
\label{eq:P_error}
    P(a,b,\varepsilon) = \begin{bmatrix}
        I_2 + O(\varepsilon^2) & O(\varepsilon)\\
        O(\varepsilon) & O(\varepsilon^2)
    \end{bmatrix}.
\end{align}
It is enough to show $\partial_{\varepsilon}P(a,b,0)$ is 0 on the diagonal blocks, for any $(a,b)$. We compute this partial. We know that $P^2 = P$ since $P$ is a projection. Take partial in $\varepsilon$ on both sides:
\begin{align}
\label{eq:equation_P_P_partial}
    \big(\partial_{\varepsilon}P(a,b,0)\big) P(a,b,0) + P(a,b,0)\big(\partial_{\varepsilon}P(a,b,0)\big) = \partial_{\varepsilon}P(a,b,0).
\end{align}
We already know $P(a,b,0) = I_k \oplus 0_{n-k}$. Say $\partial_{\varepsilon} P(a,b,0)$ has form:
\begin{align*}
    \partial_{\varepsilon} P(a,b,0) = \begin{bmatrix}
        Q & R\\
        U & V
    \end{bmatrix}.
\end{align*}
Plug into Equation~\ref{eq:equation_P_P_partial}, we get:
\begin{align*}
    \begin{bmatrix}
        Q & R\\
        0 & 0
    \end{bmatrix} + 
    \begin{bmatrix}
        Q & 0\\
        U & 0
    \end{bmatrix} =
    \begin{bmatrix}
        Q & R\\
        U & V
    \end{bmatrix},
\end{align*}
which implies both $Q$ and $V$ are 0. Hence the claimed error bound~\eqref{eq:P_error} on $P$. Then we plug in the definition $S = [Pe_1, Pe_2, (I-P)e_3, ..., (I-P)e_n]$ and obtain the claimed error bound on $S$. \\

As for $S^{-1}$, we write it as:
\begin{align*}
    S^{-1} = \left(I_n + \begin{bmatrix}
        O(\varepsilon^2) & O(\varepsilon)\\
        O(\varepsilon) & O(\varepsilon^2)
    \end{bmatrix}\right)^{-1} =: (I_n + E)^{-1} = I_n - E + O(\|E\|^2) = I_n - E + O(\varepsilon^2)
\end{align*}
by the derivative of inverse formula~\cite[p32, Equation~4.28]{Kato1995}.
\end{proof}

Then we construct similar matrix $SMS^{-1}$ using $S$. We claim $SMS^{-1}$, when viewed in a $(2\times2) \oplus ((n-2)\times(n-2))$ block structure, has error $O(\varepsilon^2)$ on the diagonal blocks and 0 on the off-diagonal blocks. We state this claim along with some additional approximations in the following proposition:

\begin{proposition}
\label{prop:block_diag_of_M}
Consider the LSE(1) setup. There exists positive constant $C>0$ such that for any $(a,b,\varepsilon)\in [1-a^*, 1+a^*]\times [-b^*, b^*]\times [0,\varepsilon^*]$, we have
\begin{align*}
    M(a,b, \varepsilon) &\sim \begin{bmatrix}
        M^{2\times 2}(a,b,\varepsilon) + \tx{Error}_1 & 0\\
        0 & M_{(n-2)\times (n-2)}(a,b,\varepsilon) + \tx{Error}_2
    \end{bmatrix},
\end{align*}
where $M^{2\times 2}$ is defined as:
\begin{align*}
    M^{2\times 2}(a,b,\varepsilon) :=  \left(I-\frac{A_{2\times2}}{a}\right)(bA_{2\times2} - aB_{2\times2})^{-1},
\end{align*}
where $A_{2\times2}$ and $B_{2\times2}$ are the top-left $2\times2$ submatrices of $A$ and $B$, respectively; while $M_{(n-2)\times(n-2)}$ refers to the bottom-right $(n-2)\times(n-2)$ submatrix of $M$. The two error terms satisfy an entry-wise second order bound:
\begin{align*}
    |(\tx{Error}_k)_{i,j}| \leq C \varepsilon^2, \quad k = 1,2
\end{align*}
for each entry $i,j$.
\end{proposition}
\begin{proof}
We already know the asymptotic structure of $S$ and $S^{-1}$ by Lemma~\ref{lemma:S_error}. Plug into $SMS^{-1}$:
\begin{align*}
    S^{-1} M S = \begin{bmatrix}
        M_{2\times 2} + O(\varepsilon^2) & O(\varepsilon)\\
        O(\varepsilon) & M_{(n-2)\times(n-2)} + O(\varepsilon^2)
    \end{bmatrix}
    =: \begin{bmatrix}
        Q & R\\
        U & V
    \end{bmatrix}.
\end{align*}
We want to show the off-diagonal blocks $R$ and $U$ are in fact 0. Recall that $S$ is defined in terms of the kernel and image of $P$ and thus $SPS^{-1} = I_2 \oplus 0_{n-2}$. Additionally, we also know $MP=PM$. We combine these equations:
\begin{align*}
    &S^{-1}MPS = (S^{-1}MS)(S^{-1}PS) = \begin{bmatrix}
        Q & 0 \\
        U & 0
    \end{bmatrix}\\
    =&\, S^{-1}PMS = (S^{-1}PS)(S^{-1}MS) = \begin{bmatrix}
        Q & R \\
        0 & 0
    \end{bmatrix}.
\end{align*}
Hence $R$ and $U$ are 0, and we conclude
\begin{align*}
    S^{-1} M S = \begin{bmatrix}
        M_{2\times 2} + O(\varepsilon^2) & 0\\
        0 & M_{(n-2)\times(n-2)} + O(\varepsilon^2)
    \end{bmatrix}.
\end{align*}

Next we show the top-left $2\times2$ block of $S^{-1}MS$ can be approximated by $M^{2\times2}$. The proof is by plain computation.
\begin{align*}
    M_{2\times 2} &:= \bigg[\left(I - \frac{A}{a}\right)(bA-aB)^{-1}\bigg]_{2\times2}\\
    &= \bigg[\left(I - \frac{A(0)}{a}\right)(bA(0) - aB(0))^{-1}\\
    &\qquad + \varepsilon \big(-\frac{A'(0)}{a}(bA(0) - aB(0))^{-1} - \left(I - \frac{A(0)}{a}\right)(bA(0) - aB(0))^{-1} (bA'(0) - aB'(0)) (bA(0)-aB(0))^{-1}\big) \\
    &\qquad + O(\varepsilon^2)\bigg]_{2\times 2}\\
    &= \left[\left(I - \frac{A(0)}{a}\right)(bA(0) - aB(0))^{-1}\right]_{2\times 2}\\
    &\qquad + \varepsilon\left[\big(-\frac{A'(0)}{a}(bA(0) - aB(0))^{-1} - \left(I - \frac{A(0)}{a}\right)(bA(0) - aB(0))^{-1} (bA'(0) - aB'(0)) (bA(0)-aB(0))^{-1}\big)\right]_{2\times 2} \\
    &\qquad + O(\varepsilon^2)\\
    &= \left(I_2 - \frac{A(0)_{2\times 2}}{a}\right)(bA(0)_{2\times 2} - aB(0)_{2\times 2})^{-1}\\
    &\qquad +\varepsilon\bigg[-\frac{A'(0)_{2\times 2}}{a} (bA(0)_{2\times 2} - aB(0)_{2\times 2})^{-1} \\
    &\qquad \qquad -\left(I_2 - \frac{A(0)_{2\times 2}}{a}\right)(bA(0)_{2\times 2} - aB(0)_{2\times 2})^{-1} (bA'(0)_{2\times 2} - aB'(0)_{2\times 2}) (bA(0)_{2\times 2}-aB(0)_{2\times 2})^{-1}\bigg]\\
    &\qquad + O(\varepsilon^2)\\
    &= M^{2\times2} + O(\varepsilon^2).
\end{align*}
\end{proof}

Now we are ready to prove the second part of the LSE setup.

\begin{theorem}[Separation of the Eigenvalues of Near Collision Line]
\label{thm:separation_eval_near_collision_line}
Consider the LSE(1) setup. There exists a $2\times2$ matrix $T(a,b,\varepsilon)$ such that the two eigenvalues of $M$ that live in $B_r(M(1,0,0))$ are precisely the eigenvalues of $T$, for any $(a,b,\varepsilon)$; while the other eigenvalues are simple and real. Consequently,
\begin{align*}
    \rho(a,b,\varepsilon) = \frac{1}{2\pi} (|\Im\,\lambda_1(T(a,b,\varepsilon))| + |\Im\,\lambda_2(T(a,b,\varepsilon))|).
\end{align*}
Furthermore, $T$ can be approximated by $M^{2\times2}$, as defined in Proposition~\ref{prop:block_diag_of_M}, entry-wise up to second order. I.e., there exists constant $C>0$ s.t. for any $(a,b,\varepsilon)$ and entry $(i,j)$, we have
\begin{align*}
    |T(a,b,\varepsilon)_{i,j} - M^{2\times2}(a,b,\varepsilon)_{i,j}| \leq C \varepsilon^2.
\end{align*}
\end{theorem}
\begin{proof}
Recall Proposition~\ref{prop:block_diag_of_M} gives a block diagonalization of $M$:
\begin{align*}
    M(a,b, \varepsilon) &\sim \begin{bmatrix}
        M^{2\times 2}(a,b,\varepsilon) + \tx{Error}_1 & 0\\
        0 & M_{(n-2)\times (n-2)}(a,b,\varepsilon) + \tx{Error}_2
    \end{bmatrix}.
\end{align*}
Let $T$ be the top-left block. We only need to show the eigenvalues of $T$ lie in $B_r(M(1,0,0))$. The rest of the statement follows from Proposition~\ref{prop:block_diag_of_M}. \\

Fix any $(a_0, b_0)$. We know we have strictly equality initially at $\varepsilon=0$:
\begin{align*}
    T(a_0,b_0,0) = M^{2\times2}(a_0,b_0,0) = M(a_0,b_0,0)_{2\times2},
\end{align*}
where $M(a_0,b_0,0)_{2\times2}$ refers to the top-left $2\times2$ submatrix of $M$. So the eigenvalues of $T(a_0,b_0,0)$ lie inside $B_r(M(1,0,0))$. If at some $\varepsilon_0$ some eigenvalue of $T(a_0,b_0,\varepsilon_0)$ leave $B_r(M(1,0,0))$, we know by continuity that there must be some $\varepsilon_1 \in [0,\varepsilon_0)$ such that one eigenvalue of $T(a_0, b_0,\varepsilon_1)$ lie on the boundary of $B_r(M(1,0,0))$, which is impossible by separation assumption of the LSE(1) setup (statement 2 of Theorem~\ref{thm:local_properties_eval_near_collision_line}).
\end{proof}

\begin{remark}
The exact expression for $T$ is given by 
\begin{align*}
    T = [SMS^{-1}]_{2\times2}.
\end{align*}
This expression is still complicated. We still do not have a good way to study $\rho$ even if we know its value only depends on $T$. The other matrix, $M^{2\times2}$, however, has a simple formula. Because $M^{2\times2}$ provides an entry-wise approximation of $T$, one would expect the eigenvalues of $M^{2\times2}$ should approximate the eigenvalues of $T$, and consequently $\rho$. We study this approximation in the following sections.
\end{remark}

\section{Width of Support and Density of $\rho$}
\label{sec:width_of_support_and_density}

We use the LSE(2) setup in this section. Recall we know $M^{2\times2}$ provides an entry-wise approximation of $T$. This section we shall upgrade it to an approximation of $\rho$ and use it to study the properties of $\rho$. The concluding theorem for this section is Theorem~\ref{thm:width_of_support_and_density}, which also proves the last part of the LSE(2) setup. Very briefly, Theorem~\ref{thm:width_of_support_and_density} states the following. For every fixed $\varepsilon \in [0,\varepsilon^*]$ such that $\rho(\cdot,\cdot,\varepsilon)$ does not varnish everywhere, the following informal statements hold,
\begin{enumerate}
    \item The boundary of $\tx{supp}(\rho(\cdot,\cdot,\varepsilon))$ is given by two curves $C_1, C_2$. 
    \item If $\rho(\cdot,\cdot,\varepsilon) > 0$ between $C_1$ ad $C_2$
    \item The distance between $C_1$ ad $C_2$ is $O(\varepsilon)$.
\end{enumerate}

\begin{figure}[H]
    \centering
    \includegraphics[width=0.5\linewidth]{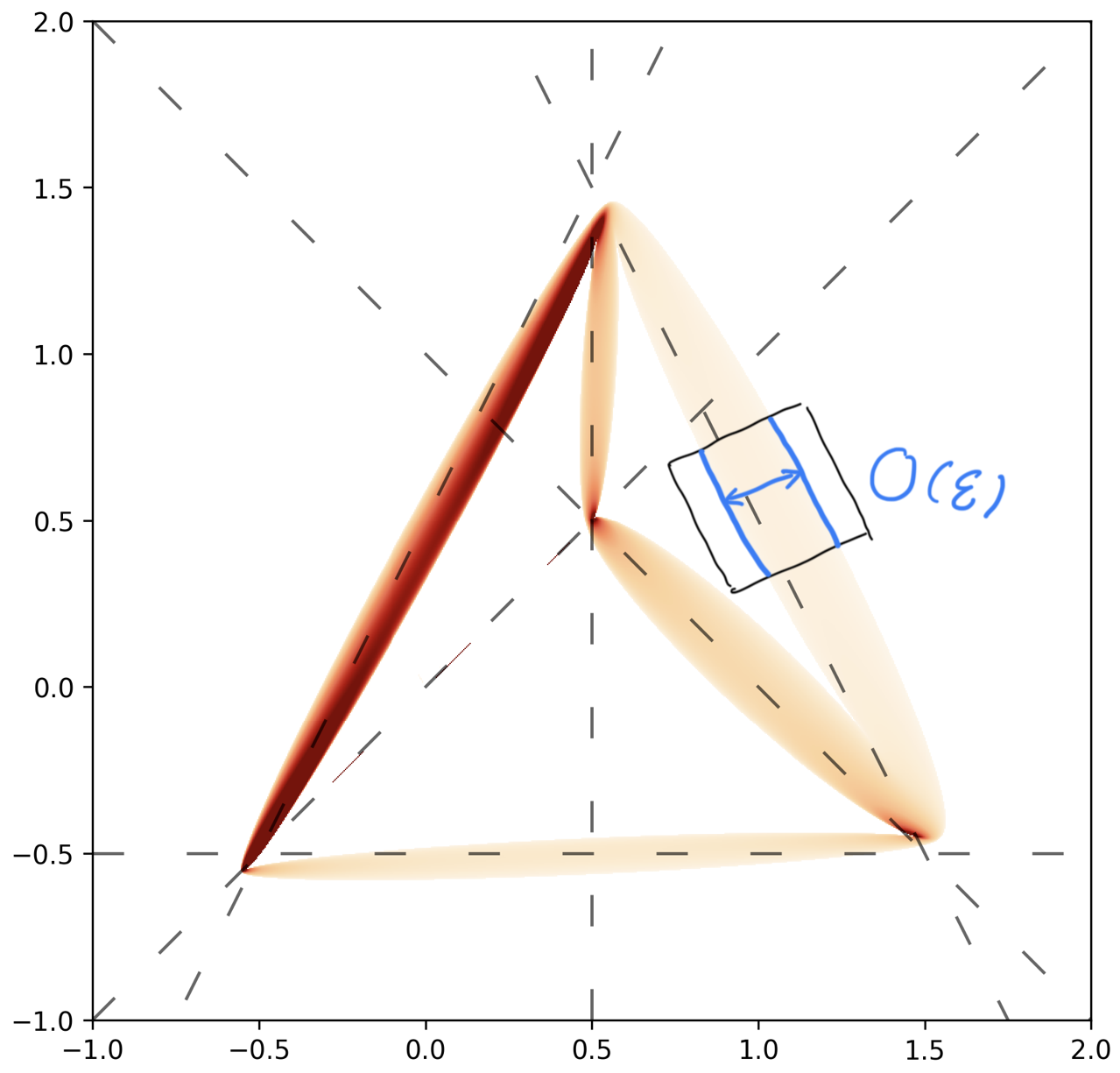}
    \caption{Illustration of part of Theorem~\ref{thm:width_of_support_and_density}. Inside an LSE rectangle (black), if $\rho$ is not identically zero, then $\partial\,\tx{supp}(\rho)$ consists of two curves (blue) separated by $O(\varepsilon)$, with $\rho>0$ between them.}
\end{figure}

Consider the discriminant $\tx{Disc}_{\lambda} \det(\lambda I - M(a,b,\varepsilon))$. For easier notation, we drop the ``$\det$'' term and simply denote:
\begin{align*}
    &\tx{Disc}(\lambda I - M): \mathbb R^3 \,\backslash \{(a,b,\varepsilon) \tx{ where } a = 0 \tx{ or } bA-aB \tx{ not invertible}\} \to \mathbb C\\
    &\tx{Disc}(\lambda I - M)(a,b,\varepsilon) := \tx{Disc}_{\lambda} \det(\lambda I - M(a,b,\varepsilon)).
\end{align*}

We define $\tx{Disc}(\lambda I - T)$ and $\tx{Disc}(\lambda I - M^{2\times2})$ in a similar way.\\

We know the following general statements on the discriminant of a real polynomial 
\begin{lemma}\cite[Section~6.2]{Barbeau1989}
Consider a polynomial $p$ with real coefficients. If all its roots are distinct, then $\tx{sign}(\tx{Disc}(p)) = (-1)^{r}$, where $r$ is the number of non-real conjugate pairs. If $p$ has repeated roots, then the discriminant is 0.
\end{lemma}

In our context, we have a real polynomial $\det(\lambda I - M(a_0,b_0,\varepsilon_0))$ with variable $\lambda$. The lemma implies that $M$ has repeated eigenvalues if and only if $\tx{Disc}(\lambda I - M)(a_0, b_0, \varepsilon_0) = 0$. Also notice that $\rho(a_0,b_0,\varepsilon_0)=0$ if and only if there is no non-real conjugate eigenvalue pairs, equivalently $\tx{Disc}(\lambda I - M)(a_0, b_0, \varepsilon_0)\geq 0$.\\

Furthermore,

\begin{proposition}
\label{prop:zeros_of_Disc_contains_boundary}
The zero set of $\tx{Disc}(\lambda I - M)$ is a superset of the boundary of $\tx{supp}(\rho)$.
\end{proposition}
\begin{proof}
Pick $(a_0,b_0,\varepsilon_0)$ that belongs to the boundary of $\tx{supp}(\rho)$. Notice $\tx{Disc}(\lambda I - M)(\cdot, \cdot,\varepsilon_0)$ is continuous in $(a,b)$. If $\tx{Disc}(\lambda I - M)(a_0,b_0,\varepsilon_0) > 0$, then there exists $U$ around $(a_0, b_0,\varepsilon_0)$ s.t. $\tx{Disc}(\lambda I - M) > 0$ in $U$, consequently $\rho>0$ everywhere in $U$, which contradicts the assumption that $(a_0, b_0,\varepsilon
_0)$ belongs to the boundary. Similarly for $\tx{Disc}(a_0,b_0,\varepsilon_0) < 0$. So it must be the case that the discriminant is 0.
\end{proof}

We specialize this proposition to the LSE(2) setup. 

\begin{corollary}
\label{coro:zeros_of_Disc_contains_boundary_LSE}
Consider the LSE(2) setup. The zero set of $\tx{Disc}(\lambda I - T)$ is a superset of $\partial\,\tx{supp}(\rho) \cap [1-a^*,1+a^*]\times[-b^*,b^*]\times[0,\varepsilon^*]$.
\end{corollary}
\begin{proof}
We know in the LSE(2) setup, the only possibly non-real eigenvalues of $M$ are the two eigenvalues of $T$, say $\lambda_1$ and $\lambda_2$, while the other $n-2$ are simple and real. Consequently,
\begin{equation*}
\begin{split}
    &(a,b,\varepsilon) \in \partial\, \tx{supp}(\rho) \cap [1-a^*, 1+a^*]\times[-b^*,b^*]\times[0,\varepsilon^*]\\
    \implies&\tx{Disc}(\lambda I - M)(a, b, \varepsilon) = 0 \\
    \iff& \lambda_1(a, b, \varepsilon) = \lambda_2(a, b, \varepsilon)\\
    \iff& \tx{the two eigenvalues of } T := M^{2\times 2} + \tx{Error} \tx{ are equal}\\
    \iff& \tx{Disc}(\lambda I - T)(a,b,\varepsilon) = 0.
\end{split}
\end{equation*}
\end{proof}

\subsection{Zero Set of $\tx{Disc} (\lambda I - T)$}

We continue on the LSE(2) setup and Corollary~\ref{coro:zeros_of_Disc_contains_boundary_LSE}. We want to understand the zero set of $\tx{Disc}(\lambda I - T)$ and also $\tx{Disc}(\lambda I - M^{2\times2})$ since $M^{2\times2}$ approximates $T$. Notice the matrices are $2\times2$ and we may expand the respective discriminants directly:
\begin{align*}
    &\tx{Disc}(\lambda I - T) = (T_{1,1}-T_{2,2})^2 + 4T_{1,2}T_{2,1}\\
    &\tx{Disc}(\lambda I - M^{2\times2}) = (M^{2\times2}_{1,1}-M^{2\times2}_{2,2})^2 + 4M^{2\times2}_{1,2}M^{2\times2}_{2,1}.
\end{align*}

\begin{lemma}
\label{lemma:M_entry_Disc_bounds}
Consider the LSE(2) setup. There exists $C>0$ such that for any $(a,b,\varepsilon)$, the following holds:
\begin{enumerate}
    \item $|M^{2\times2}_{1,1} - M^{2\times2}_{2,2}| \leq C(|a-1|+\varepsilon)$.
    \item $|M^{2\times2}_{1,2}|, |M^{2\times2}_{2,1}| \leq C\varepsilon$.
    \item $|\tx{Disc}(\lambda I - T)(a, b, \varepsilon) - \tx{Disc}(\lambda I - M^{2\times2})(a,b,\varepsilon)| \leq C \varepsilon^2(|a-1| + \varepsilon)$.
    \item If $(a,b,\varepsilon) \in \partial\,\tx{supp}(\rho)$, then it must be the case that $|\tx{Disc}(\lambda I - M^{2\times2})(a,b,\varepsilon)| \leq C \varepsilon^2 (|a-1|+\varepsilon)$.
\end{enumerate}
\end{lemma}
\begin{proof}
First we prove statements (1) and (2). Denote $f:=M^{2\times2}_{1,1} - M^{2\times2}_{2,2}$, we have:
\begin{align*}
    |f(a,b,\varepsilon)| &\leq \underbrace{|f(1,b,0)|}_{=0} + |a-1| \max_{[1-a^*,1+a^*]} |\partial_{x} f(x, b, \varepsilon)| + \varepsilon \max_{[0,\varepsilon^*]} |\partial_z f(a,b,z)|.
\end{align*}
For $|M^{2\times2}_{1,2}|$:
\begin{align*}
    |M^{2\times2}_{1,2}(a,b,\varepsilon)| &\leq \underbrace{|M^{2\times2}_{1,2}(a,b,0)|}_{=0} + \varepsilon \max_{[0,\varepsilon^*]} |\partial_{z}M^{2\times2}_{1,2}(a,b_0,z)|.
\end{align*}
Similarly for $|M^{2\times2}_{2,1}|$. We choose constant $C_1$ to be max of modulus of all partials involved.\\

Then we prove statement (3). Take the difference of the two discriminants:
\begin{align*}
    &\tx{Disc}(\lambda I - T) - \tx{Disc}(\lambda I - M^{2\times2})\\
    =&\, 2(M^{2\times2}_{1,1}-M^{2\times2}_{2,2})
    \big((T_{1,1}-M^{2\times2}_{1,1})-(T_{2,2}-M^{2\times2}_{2,2})\big) \\
    &+\big((T_{1,1}-M^{2\times2}_{1,1})-(T_{2,2}-M^{2\times2}_{2,2})\big)^2\\
    &+4M^{2\times2}_{1,2}(T_{2,1}-M^{2\times2}_{2,1})
    +4M^{2\times2}_{2,1}(T_{1,2}-M^{2\times2}_{1,2}) +4(T_{1,2}-M^{2\times2}_{1,2})(T_{2,1}-M^{2\times2}_{2,1}).
\end{align*}
Every term $|T_{i,j} - M^{2\times2}_{i,j}|$ is bounded by $C_1\varepsilon^2$ by Theorem~\ref{thm:separation_eval_near_collision_line}, so
\begin{align}
\label{ineq:Disc_bounds_intermediate}
    |\tx{Disc}(\lambda I - T) - \tx{Disc}(\lambda I - M^{2\times2})|
    \leq C\varepsilon^2(2|M^{2\times2}_{1,1}-M^{2\times2}_{2,2}| + 4|M^{2\times2}_{1,2}|+4|M^{2\times2}_{2,1}|) + 8C\varepsilon^4.
\end{align}
We know all of $|M^{2\times2}_{1,1}-M^{2\times2}_{2,2}|$, $|M^{2\times2}_{1,2}|$, $|M^{2\times2}_{2,1}|$ can be bounded by $C_1(|a-1| + \varepsilon)$ by statements (1) and (2). Plugging this bound into Inequality~\ref{ineq:Disc_bounds_intermediate} gives
\begin{align*}
    |\tx{Disc}(\lambda I - T)(a, b, \varepsilon) - \tx{Disc}(\lambda I - M^{2\times2})(a,b,\varepsilon)| \leq C_2 \varepsilon^2(|a-1| + \varepsilon),
\end{align*}
for some constant $C_2$. Setting $C = \max\{C_1, C_2\}$ so that all three statements hold at the same time.\\

Lastly we prove statement (4). Say $(a,b,\varepsilon) \in \partial\,\tx{supp}(\rho)$. By Corollary~\ref{coro:zeros_of_Disc_contains_boundary_LSE}, we have $\tx{Disc}(\lambda I - T)(a,b_0,\varepsilon) = 0$; plugging into statement (3) gives the desired bound on $\tx{Disc}(\lambda I - M^{2\times2})$.
\end{proof}

\begin{proposition}[Width of Support]
\label{prop:v_bound_by_eps_at_zero}
Consider the LSE(2) setup. There exists $C>0$ such that for any $b_0 \in [-b^*,b^*]$ and any $(a,b_0,\varepsilon) \in \partial\,\tx{supp}(\rho(\cdot,b_0,\cdot))$, it must be the case that $|a-1| \leq C \varepsilon$.
\end{proposition}
\begin{figure}[H]
    \centering
    \includegraphics[width=0.35\linewidth]{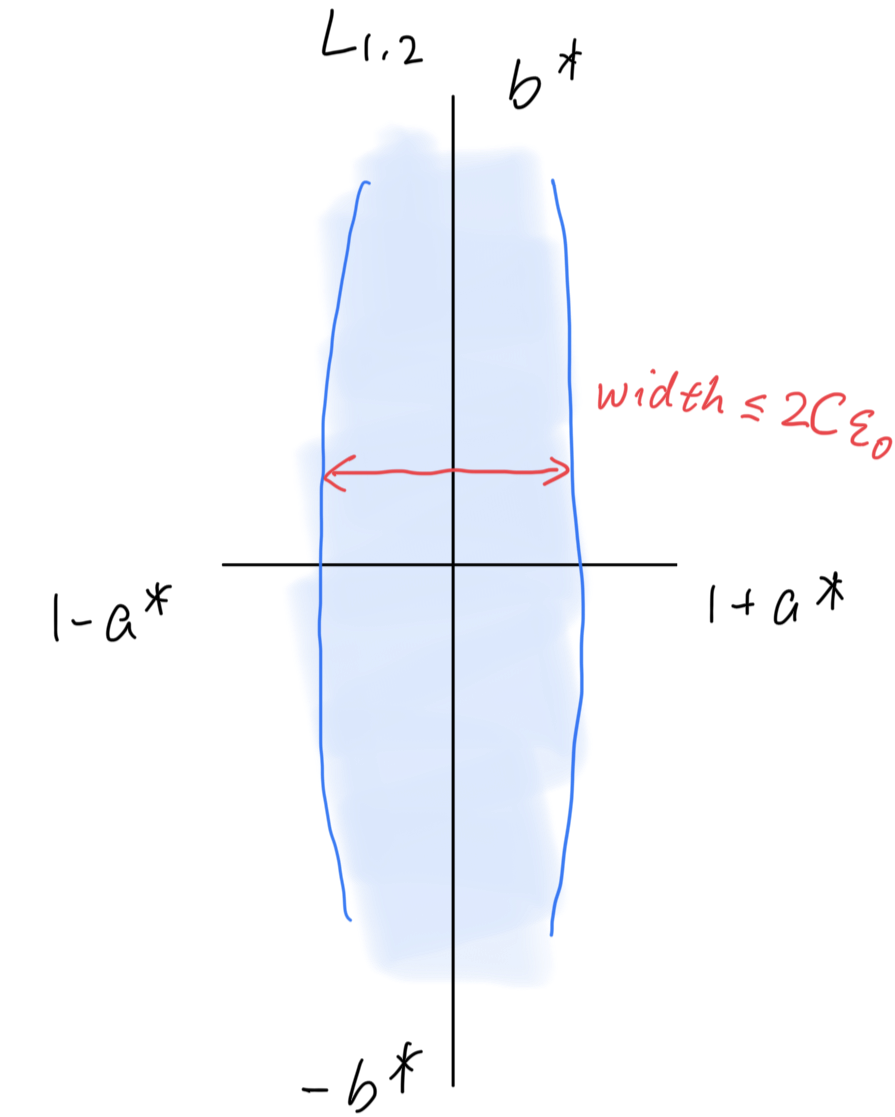}
    \caption{The width of support is first order in $\varepsilon$. We fix some $\varepsilon_0$ and draw this illustration.}
    \label{fig:width_of_support}
\end{figure}
\begin{proof}
Fix some $b_0$ and pick $(a,b_0,\varepsilon) \in \partial\,\tx{supp}(\rho(\cdot,b_0,\cdot))$. On the one hand, we know by Lemma~\ref{lemma:M_entry_Disc_bounds} (4) that
\begin{align*}
    &|\tx{Disc}(\lambda I - M^{2\times2})| \leq C_1 \varepsilon^2 (|a-1|+\varepsilon)\\
    \iff& |(M^{2\times2}_{1,1}-M^{2\times2}_{2,2})^2 + 4M^{2\times2}_{1,2}M^{2\times2}_{2,1}| \leq C_1 \varepsilon^2 (|a-1|+\varepsilon)\\
    \implies& (M^{2\times2}_{1,1} - M^{2\times2}_{2,2})^2 \leq 4|M^{2\times2}_{1,2} M^{2\times2}_{2,1}| + C_1 \varepsilon^2(|a-1|+\varepsilon)\\
    \iff& |M^{2\times2}_{1,1} - M^{2\times2}_{2,2}| \leq \sqrt{4|M^{2\times2}_{1,2} M^{2\times2}_{2,1}| + C_1 \varepsilon^2(|a-1|+\varepsilon)}.
\end{align*}
We know by Lemma~\ref{lemma:M_entry_Disc_bounds}(2) that there exists $C_2>0$ such that both $|M^{2\times2}_{1,2}|$ and $|M^{2\times2}_{2,1}|$ are upper bounded by $C_2\varepsilon$. Plug into the inequality above:
\begin{align}
\label{ineq:diff_M(eps)_upper_bound}
    |M^{2\times2}_{1,1} - M^{2\times2}_{2,2}| &\leq \sqrt{4 C_2^2 \varepsilon^2 + C_1 \varepsilon^2(|a-1|+\varepsilon)} \leq \varepsilon\sqrt{4C_2^2 + C_1(a^*+1 +\varepsilon^*)}
\end{align}
which holds for all $(a,\varepsilon)$. On the other hand, Fundamental Theorem of Calculus gives:
\begin{align*}
    &\left|M^{2\times2}_{1,1}(a, b_0,\varepsilon) - M^{2\times2}_{2,2}(a, b_0, \varepsilon) - \big(M^{2\times2}_{1,1}(a,b_0,0) - M^{2\times2}_{2,2}(a, b_0, 0)\big)\right| \leq C_3 \varepsilon\\
    \implies& \left|M^{2\times2}_{1,1}(a,b_0,0) - M^{2\times2}_{2,2}(a,b_0,0)\right| \leq C_3 \varepsilon + \left| \big(M^{2\times2}_{1,1}(a,b_0,\varepsilon) - M^{2\times2}_{2,2}(a,b_0,\varepsilon)\big)\right|.
\end{align*}
Combine with Inequality~\eqref{ineq:diff_M(eps)_upper_bound}:
\begin{align}
\label{ineq:diff_M(0)_upper_bound}
    |M^{2\times2}_{1,1}(a,b_0,0) - M^{2\times2}_{2,2}(a,b_0,0)| &\leq \varepsilon\left(C_3 + \sqrt{4C_1^2 + C_2(a^*+1+ \varepsilon^*)}\right) =: \varepsilon C_4.
\end{align}
Notice $M^{2\times2}(a,b_0,0)$ is simply a $2\times2$ diagonal matrix and we may compute the difference of its diagonal terms directly:
\begin{align*}
    M^{2\times2}_{1,1}(a,b_0,0) - M^{2\times2}_{2,2}(a,b_0,0) &= \frac{b_0(D_{2,2}-D_{1,1}) - a(E_{2,2}-E_{1,1})+D_{1,1}E_{2,2}-D_{2,2}E_{1,1}}{(b_0D_{1,1}-aE_{1,1})(b_0D_{2,2}-aE_{2,2})}.
\end{align*}
We know the denominator is continuous in $a$ and non-zero everywhere on $[1-a^*, 1+a^*]$ (otherwise $M(a,b_0,0)$ would not be well-defined). So there exists some $C_5$ such that the absolute value of the denominator is uniformly lower bounded by $C_5\varepsilon$. As for the numerator, we know it is affine in $a$ and equal to 0 at $a=0$ (the two entries $M^{2\times2}_{1,1}(0,b_0,0)$ and $M^{2\times2}_{2,2}(0,b_0,0)$ are the eigenvalues of $M^{2\times2}(0,b_0,0)$ and are equal on the collision line). So the numerator is linear in $a$, say it equals $C_6 a$. Plug in:
\begin{align*}
    \frac{C_6|a|}{C_5}  &\leq |M^{2\times2}_{1,1}(a,b_0,0) - M^{2\times2}_{2,2}(a,b_0,0)| \leq \varepsilon C_4,
\end{align*}
which says $|a|$ is upper bounded by some $C\varepsilon$.
\end{proof}

\begin{proposition}[Shape of Support, Non-Degenerate Case]
\label{prop:Disc_at_most_two_zeros}
Consider the LSE(2) setup. There exists $0<\varepsilon^*_{cvx}\leq \varepsilon^*$ such that for every $\varepsilon_0\in [0,\varepsilon^*_{cvx}]$ and $b_0 \in [-b^*, b^*]$, $\tx{Disc}(\lambda I - T(a, b_0, \varepsilon_0))$ is strictly convex in $a$ and has at most two zeros on $[1-a^*, 1+a^*]$. If there are indeed two distinct $a$'s s.t. $\tx{Disc}(\lambda I - T) = 0$, say $a_-$ and $a_+$, then $(a_-, b_0), (a_+, b_0) \in \partial\,\tx{supp}(\rho(\cdot, \cdot, \varepsilon_0))$ and $\rho$ is positive on the segment between $a_-$ and $a_+$, zero outside. 
\end{proposition}
\begin{proof}
We use the fact that a strictly convex function has at most two zeros. This is half of the statement; the other half is, the strictly convex function is strictly negative between the two zeros, and strictly positive outside.\\

Fix some $b_0 \in [-b^*, b^*]$. First we prove $\tx{Disc}(\lambda I - T(a,b_0,0))$, is strictly convex at in $a$ at $1$, and then we pass the strict convexity to a neighborhood around $(1,b_0,0)$. Recall we have proved an error bound on the discriminants of $T$ and $M^{2\times2}$ (Lemma~\ref{lemma:M_entry_Disc_bounds}(3)):
\begin{align*}
    |\tx{Disc}(\lambda I - T(a,b_0,\varepsilon)) - \tx{Disc}(\lambda I - M^{2\times2}(a,b_0,\varepsilon))| \leq C_2 \varepsilon^2(|a| + \varepsilon).
\end{align*}
In particular, we have strict identity at $\varepsilon=0$:
\begin{align*}
    \tx{Disc}(\lambda I - T)(a,b_0,0) &= \tx{Disc}(\lambda I - M^{2\times2})(a,b_0,0)\\
    &= (M^{2\times2}_{1,1}(a,b_0,0) - M^{2\times2}_{2,2}(a,b_0,0))^2 + 4 M^{2\times2}_{1,2}(a,b_0,0) M^{2\times2}_{2,1}(a,b_0,0)\\
    &= (M^{2\times2}_{1,1}(a,b_0,0) - M^{2\times2}_{2,2}(a,b_0,0))^2
\end{align*}
\begin{align*}
    \implies \partial_a^2\tx{Disc}(\lambda I - T)(a,b_0,0) &= \,2\big(\partial_aM^{2\times2}_{1,1}(a,b_0,0) - \partial_aM^{2\times2}_{2,2}(a,b_0,0)\big)^2\\
    &\quad + 2(M^{2\times2}_{1,1}(a,b_0,0) - M^{2\times2}_{2,2}(a,b_0,0))(\partial_a^2M^{2\times2}_{1,1}(a,b_0,0) - \partial_a^2M^{2\times2}_{2,2}(a,b_0,0)).
\end{align*}
Now we further consider $a=1$. At $a=1,\,\varepsilon=0$, we know the diagonal entries of $M^{2\times2}$ are equal, hence:
\begin{align*}
    \partial_a^2\tx{Disc}(\lambda I - T)(1,b_0,0) &= \,2\big(\partial_aM^{2\times2}_{1,1}(1,b_0,0) - \partial_aM^{2\times2}_{2,2}(1,b_0,0)\big)^2.
\end{align*}
We compute this:
\begin{align*}
    M^{2\times2}_{1,1}(a,b_0,0) - M^{2\times2}_{2,2}(a,b_0,0) = \frac{b_0(D_{2,2} - D_{1,1}) - a(E_{2,2}-E_{1,1})+D_{1,1}E_{2,2}-D_{2,2}E_{1,1}}{(b_0 D_{1,1}- aE_{1,1})(b_0 D_{2,2}-aE_{2,2})}.
\end{align*}
The numerator varnishes at $a=1$ by the assumption that $(1,b_0)\in L_{1,2}$. So we have:
\begin{align*}
    \big(\partial_a M^{2\times2}_{1,1}(1,b_0,0) - \partial_a M^{2\times2}_{2,2}(1,b_0,0)\big)^2 
    &=\bigg(\left.\partial_a \frac{b_0(D_{2,2} - D_{1,1}) - a(E_{2,2}-E_{1,1})+D_{1,1}E_{2,2}-D_{2,2}E_{1,1}}{(b_0D_{1,1}- aE_{1,1})(b_0D_{2,2}-aE_{2,2})}\right|_{a=1}\bigg)^2\\
    &=\bigg(\left.\frac{(\partial_a \tx{numerator}) \cdot\tx{denominator}}{\tx{denominator}^2}\right|_{a=1}\bigg)^2\\
    &= \bigg(\left.\frac{(E_{2,2}-E_{1,1}) \cdot\tx{denominator}}{\tx{denominator}^2}\right|_{a=1}\bigg)^2.
\end{align*}
The denominator is non-zero by the guarantee that $M$ is well-defined in the domain. And $E_{2,2}\neq E_{1,1}$ because $(1,E_{1,1}),\,(1,E_{2,2})$ are the two points used to define $L_{1,2}$ and are guaranteed to be distinct. So we may conclude that $\tx{Disc}(\lambda I - T)$ is strictly convex at $(1,b_0,0)$. By analyticity of the discriminant, it should remain strictly convex in a neighborhood around $(1,b_0,0)$, say $[1-a^{**},1+a^{**}]\times [-b^{**}, b^{**}]\times [0,\varepsilon_{cvx}^*]$. We wish to shrink the range of $\varepsilon$ so that the range of $a$ may grow to the full $[-a^*, a^*]$.\\

Recall Proposition~\ref{prop:v_bound_by_eps_at_zero}: for any zero $(a,b,\varepsilon)$ of $\tx{Disc}(\lambda I - T)$ in the whole rectangle $[1-a^*,1+a^*]\times[-b^*, b^*]\times[0,\varepsilon^*]$, it must be the case that $|a|\leq C\varepsilon$, where the constant $C$ is uniform. In particular we pick some $\varepsilon^*_{cvx} \in (0, \min(\varepsilon_{cvx}^*, v^{**}/C)]$. For every zero $(b,\varepsilon)$ where $\varepsilon\in [0,\varepsilon^*_{cvx}]$, we know:
\begin{align*}
    |a| \leq C\varepsilon \leq C\varepsilon^*_{cvx}\leq v^{**}.
\end{align*}
In other words, for any $\varepsilon\in [0,\varepsilon^*_{cvx}]$, we are guaranteed that no zero can appear outside $[-v^{**}, v^{**}]$. And we have concluded that there are at most two zeros inside $[-v^{**}, v^{**}]$. So there are at most two zeros in the whole interval $[-v^*, v^*]$.\\

Lastly, say there are indeed two distinct zeros. By strict convexity, $\tx{Disc}(\lambda I - T)$ is negative between the two zeros and positive outside. We pass this positivity to density $\rho$. Recall for throughout the whole rectangle $[1-a^*,1+a^*]\times[-b^*,b^*]\times [0,\varepsilon^*]$, we have:
\begin{align*}
    \rho(a,b_0,\varepsilon_0) = \frac{1}{2\pi}(|\Im\lambda_1(a,b_0,\varepsilon_0)| + |\Im\lambda_2(a,b_0,\varepsilon_0)|),
\end{align*}
which are the two eigenvalues enclosed by the ball $B_r(M(1,0,0)_{1,1})$ and gthe other eigenvalues are guaranteed to be real. Whereas the discriminant of $T$ has formula:
\begin{align*}
    \tx{Disc}(\lambda I - T)(a,b_0,\varepsilon_0) = (\lambda_1(a,b_0,\varepsilon_0) - \lambda_2(a,b_0,\varepsilon_0))^2
\end{align*}
which is negative if and only if $\lambda_1$ and $\lambda_2$ form non-trivial conjugate pairs. Given that $\tx{Disc}(\lambda I - T)$ is negative between the two zeros $a_-$ and $a_+$ and positive outside, we know $\lambda_1$ and $\lambda_2$ form non-trivial conjugate pairs in between but not outside. Hence $\rho>0$ between $a_-$ and $a_+$ and equal to 0 outside.
\end{proof}

\begin{corollary}[Shape of Support, Degenerate Case]
\label{coro:Disc_one_zero}
Consider the LSE(2) setup and $\varepsilon_{cvx}^*$ as given by Proposition~\ref{prop:Disc_at_most_two_zeros}. Shrink $\varepsilon_{cvx}^*$ if necessary s.t. $\tx{Disc}(\lambda I - T)(\cdot, b, \varepsilon)$ is strictly positive at $1-a^*$ and $1+a^*$ for every $(b,\varepsilon)\in [-b^*, b^*]\times[0,\varepsilon_{cvx}^*]$. Fix some $(b_0,\varepsilon_0)$. If $\tx{Disc}(\lambda I - T)(\cdot, b_0, \varepsilon_0)$ has at most one zero on $[1-a^*, 1+a^*]$, then it is non-negative everywhere and $\rho(\cdot, b_0, \varepsilon_0) = 0$ throughout $[1-a^*, 1+a^*]$.
\end{corollary}
\begin{proof}
It is possible to shrink $\varepsilon_{cvx}^*$ so that $\tx{Disc}(\lambda I - T)(\cdot, b, \varepsilon)$ is strictly positive at $1-a^*$ and $1+a^*$ because, both $\tx{Disc}(\lambda I - T)(1-a^*, b, 0)$ and $\tx{Disc}(\lambda I - T)(1+a^*, b, 0)$ are positive and the discriminant is continuous in $\varepsilon$.\\

Fix any $(b_0, \varepsilon_0)\in [-b^*, b^*]\times[0,\varepsilon_{cvx}^*]$, we know $\tx{Disc}(\lambda I - T)(\cdot,b_0, \varepsilon_0)$ is strictly convex on $[1-a^*, 1+a^*]$. Because of the strict convexity, the discriminant is negative only between its (at most two) zeros. If the discriminant has no zero, then it is non-negative throughout $[1-a^*, 1+a^*]$. Consequently there is no non-real eigenvalue pairs and $\rho(\cdot, b_0, \varepsilon_0) = 0$ everywhere on $[1-a^*, 1+a^*]$.
\end{proof}

\subsection{Width of Support and Density of $\rho$}

We add one more proposition on the density of $\rho$, and then we are ready to conclude this section.

\begin{proposition}
\label{prop:density_first_order_bound}
Consider the LSE(2) setup. There exists continuous function $C:[1-a^*,1+a^*]\times [-b^*,b^*] \to \mathbb R_{\geq 0}$ such that $\rho(a,b,\varepsilon) < C(a,b)\,\varepsilon$ everywhere in $[1-a^*, 1+a^*]\times[-b^*, b^*]\times [0,\varepsilon^*]$.
\end{proposition}

The proof of this proposition requires a lemma by Bauer and Fike: 
\begin{lemma}[Bauer--Fike at $l_2$ norm, {\cite[Theorem IIIa]{BauerFike1960}}]
Let $A$ be a diagonalizable matrix with diagonalization $A = V D V^{-1}$. Consider $A + \Delta$. For every eigenvalue $\mu$ of $A+\Delta$, there is an eigenvalue $\lambda$ of $A$ s.t.
\begin{align*}
    |\lambda - \mu| \leq \|V\| \|V^{-1}\| \|\Delta\|.
\end{align*}
\end{lemma}

\begin{proof}[Proof for Proposition~\ref{prop:density_first_order_bound}]

We begin with a pointwise argument, and then upgrade it to a uniform statement by compactness of $\tx{supp}(\rho)$. Pick any $(a,b,\varepsilon)\in [1-a^*,1+a^*]\times [-b^*, b^*]\times[0,\varepsilon^*]$. We apply Bauer-Fike with $A = M(a,b,0)$ and $A+\Delta = M(a,b,\varepsilon)$. In particular, we know $M(a,b,0)$ is already diagonal, so the diagonalization matrix $V$ in Bauer-Fike can be chosen to be $I$. The theorem states that for any $\lambda_i(a,b,\varepsilon)$, there exists $\lambda_j(a,b,0)$ such that:
\begin{align*}
    |\lambda_i(a,b,\varepsilon) - \lambda_j(a,b,0)| &\leq \|I \| \|I^{-1} \|M(a,b,\varepsilon) - M(a,b,0)\\
    &= \|M(a,b,\varepsilon) - M(a,b,0)\|.
\end{align*}
Recall that $M(a,b,\varepsilon)$ have at most two non-real eigenvalues, which live in $B_r(M(1,0,0)_{1,1})$ and we denote by $\lambda_1(a,b,\varepsilon)$ and $\lambda_2(a,b,\varepsilon)$. Bauer-Fike states that there exists eigenvalues $\lambda_{j_1}(a,b,0)$ and $\lambda_{j_2}(a,b,0)$ (which are real) that are close to $\lambda_1(a,b,\varepsilon)$ and $\lambda_2(a,b,\varepsilon)$. Hence
\begin{align*}
    \rho(a,b,\varepsilon) &= \frac{1}{2\pi} \left(|\Im \lambda_1(a,b,\varepsilon)| + |\Im \lambda_2(a,b,\varepsilon)|\right)\\
    &\leq \frac{1}{2\pi}(|\lambda_1(a,b,\varepsilon) - \lambda_{j_1}(a,b,0)| + |\lambda_2(a,b,\varepsilon) - \lambda_{j_2}(a,b,0)|)\\
    &= \frac{1}{\pi} \|M(a,b,\varepsilon) - M(a,b,0)\|.
\end{align*}
The last line is first order in $\varepsilon$ because $M$ is analytic.
\end{proof}

Now we may are ready to conclude the properties on the size of $\tx{supp}(\rho)$ and the density. The following theorem also provides the final statement of the LSE setup.

\begin{theorem}[Width of Support and Density]
\label{thm:width_of_support_and_density}
Consider the LSE(2) setup. The following holds inside the rectangle $[1-a^*, 1+a^*]\times[-b^*, b^*]\times [0,\varepsilon^*]$:
\begin{enumerate}
    \item (Prop~\ref{prop:v_bound_by_eps_at_zero}) There exists $C_1>0$ such that for every $b\in[-b^*,b^*]$ and every $(a,\varepsilon) \in \partial\,\tx{supp}(\rho(\cdot, b, \cdot))$, we have $|a|\leq C_1 \varepsilon$.
    \item (Prop~\ref{prop:density_first_order_bound}) There exists $C_2>0$ such that for every $(a,b,\varepsilon) \in [1-a^*,1+a^*]\times[-b^*,b^*]\times[0,\varepsilon^*]$, we have $\rho(a, b, \varepsilon) \leq C_2 \varepsilon$.
\end{enumerate}
Furthermore, there exists $0<\varepsilon_{cvx}^*\leq \varepsilon^*$ s.t. the following holds inside the smaller rectangle $[1-a^*, 1+a^*]\times[-b^*, b^*]\times [0,\varepsilon_{cvx}^*]$:
\begin{enumerate}
\setcounter{enumi}{2}
    \item (Prop~\ref{prop:Disc_at_most_two_zeros}) For any $b_0\in[-b^*, b^*]$ and $\varepsilon_0 \in  [0,\varepsilon_{cvx}^*]$, $\partial\,\tx{supp}(\rho(\cdot, b_0, \varepsilon_0))$ intersect $[1-a^*, 1+a^*]$ at two point.
    \item (Prop~\ref{prop:Disc_at_most_two_zeros}) If $\partial\,\tx{supp}(\rho(\cdot, b_0, \varepsilon_0))$ does intersect $[1-a^*, 1+a^*]$ at two distinct points, say $a_-$ and $a_+$, then we know $\rho$ is strictly positive on $(a_-, a_+)$ and zero outside.
    \item (Coro~\ref{coro:Disc_one_zero}) Otherwise $\rho(\cdot, b_0, \varepsilon_0) = 0$ everywhere on $[1-a^*, 1+a^*]$.
\end{enumerate}
\end{theorem}

\section{$2\times2$ Approximation of $\rho$}
\label{sec:22_approx_rho}

We have proved the full LSE setup in the sections above. Recall by the LSE setup, we may write:
\begin{align*}
    \rho\; \overset{\tx{Thm~\ref{thm:separation_eval_near_collision_line}}}{=}\; \frac{1}{2\pi}(|\Im \lambda_1(T)| + |\Im \lambda_2(T)|)
\end{align*}
for some $2\times2$ matrix $T$. Furthermore, we have another $2\times2$ matrix $M^{2\times2}$ with simple formula that approximates $T$ entrywise. We shall upgrade this approximation of matrix to an approximation of density:
\begin{align*}
    \rho\; \overset{\tx{LSE}}{=}\; \frac{1}{2\pi}(|\Im \lambda_1(T)| + |\Im \lambda_2(T)|) \overset{\tx{WTS}}{\approx} \frac{1}{2\pi}(|\Im \lambda_1(M^{2\times2})| + |\Im \lambda_2(M^{2\times2})|.
\end{align*}

And then upgrade this local approximation to a global version. More concretely, we conclude with the following theorem:

\begin{theorem}[$2\times2$ Approximation of Density]
\label{thm:22_approx}
There exists $\varepsilon^*>0$ and $C>0$ such that for any $(a,b,\varepsilon) \in (\mathbb R^3 / \Omega_r) \times [0,\varepsilon^*]$, we have:
\begin{align*}
    |\rho(a,b,\varepsilon) - \rho_{approx}(a,b,\varepsilon)| \leq C \varepsilon^{3/2}
\end{align*}
where the approximating density is defined as:
\begin{align*}
    \rho_{approx} &:= \sum_{i<j} \rho_{i,j}\\
    &:= \sum_{i<j}\frac{1}{2\pi}\big(|\Im\, \lambda_1(M^{\{i,j\}}(a,b,\varepsilon))| + |\Im\, \lambda_2(M^{\{i,j\}}(a,b,\varepsilon))|\big),
\end{align*}
and 
\begin{align*}
    M^{\{i,j\}}(a,b,\varepsilon) := \left(I-\frac{A_{\{i,j\}}(\varepsilon)}{a}\right)(bA_{\{i,j\}}(\varepsilon) - aB_{\{i,j\}}(\varepsilon))^{-1}
\end{align*}
where $A_{\{i,j\}}$ and $B_{\{i,j\}}$ refer to refer to the $2\times2$ principal submatrices of $A$ and $B$ corresponding to rows/columns $i$ and $j$, and $\Omega_r$ is defined as below.
\end{theorem}

\begin{definition}[Complement of Region of Approximation]
\label{def:Omega}
Pick a set of radii 
\begin{align*}
    r = \{r_{i,j,k,l} > 0 \,|\, i<j, k<l \in [1,n]\},
\end{align*}
and denote $\Omega_r$ as the following set in $\mathbb R^3$:
\begin{align*}
    \Omega_r := \Sigma &\cup \{(a,b,\varepsilon) \,|\, a = 0 \tx{ or } \det(bA_{\{i,j\}} - aB_{\{i,j\}}) = 0 \tx{ for some } 1\leq i<j\leq n\}\\
    &\cup \bigcup_{i<j,\,k<l} B_{r_{i,j,k,l}}(L_{i,j} \cap L_{k,l}).
\end{align*}
where $A_{\{i,j\}}$ and $B_{\{i,j\}}$ refer to refer to the $2\times2$ principal submatrices of $A$ and $B$ corresponding to rows/columns $i$ and $j$, and $B_r(L_{i,j} \cap L_{k,l})$ denotes the open ball centered at the point of intersection of $L_{i,j}$ and $L_{k,l}$ with radius $r_{i,j,k,l}$.
\end{definition}

\begin{figure}
    \centering
    
    \begin{subfigure}[b]{0.48\textwidth}
        \centering
        \includegraphics[width=\linewidth]{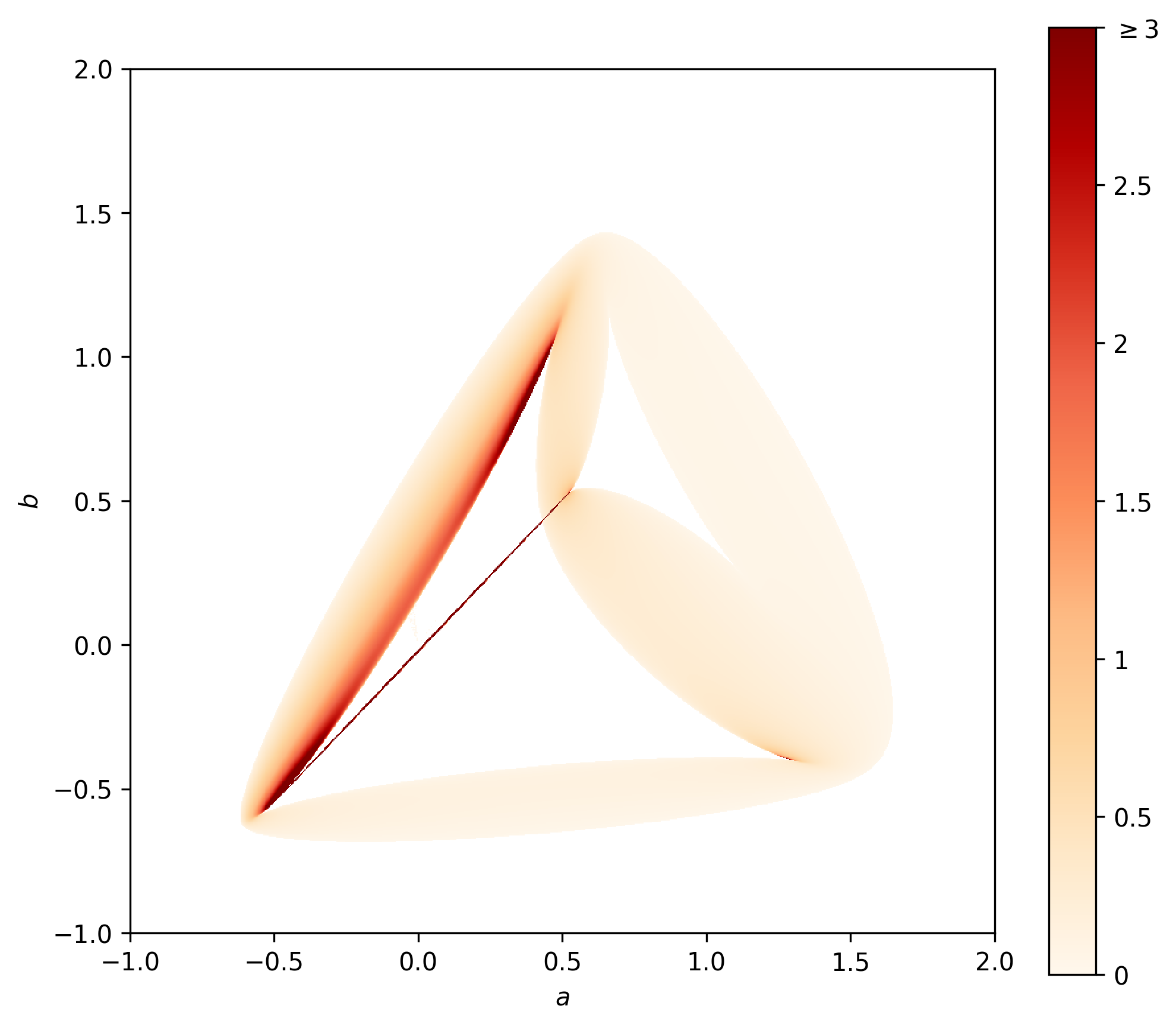}
        \caption{The true density $\rho$.}
    \end{subfigure}
    \hfill
    \begin{subfigure}[b]{0.48\textwidth}
        \centering
        \includegraphics[width=\linewidth]{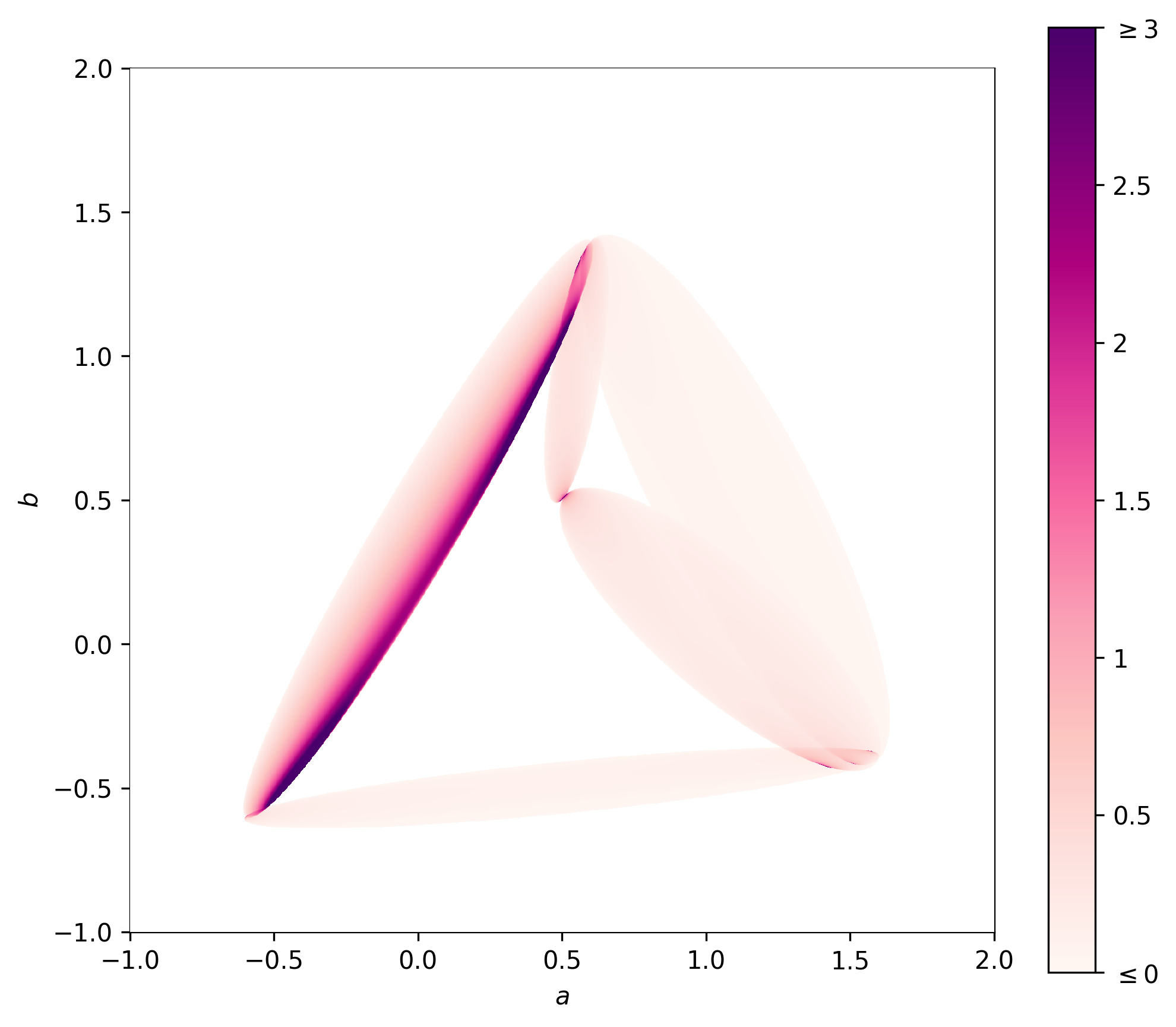}
        \caption{The approximating density $\rho_{approx}$.}
    \end{subfigure}

    \begin{subfigure}[b]{0.53\textwidth}
        \centering
        \includegraphics[width=\linewidth]{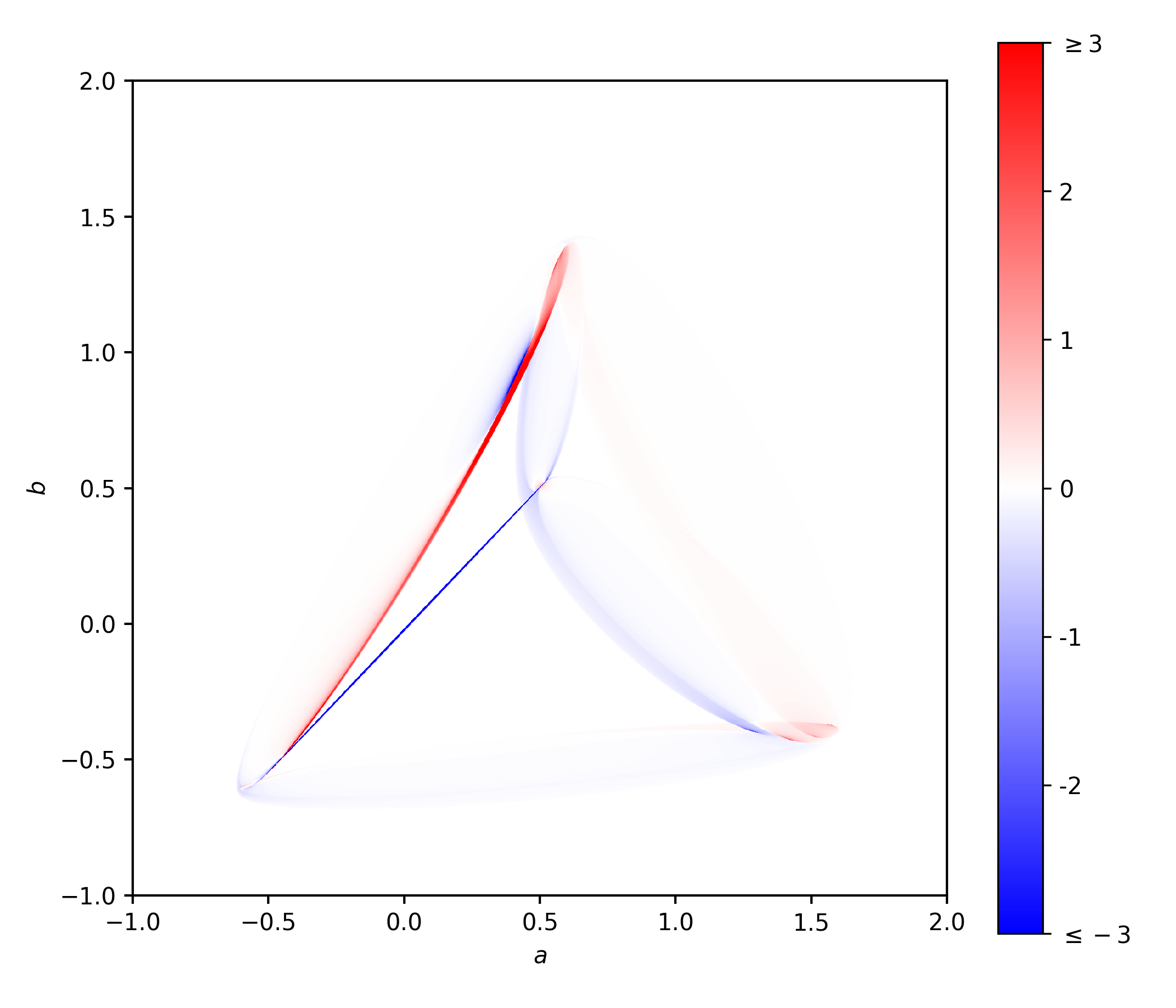}
        \caption{Error $\rho - \rho_{approx}$.}
    \end{subfigure}
    \hfill
    \begin{subfigure}[b]{0.45\textwidth}
        \centering
        \includegraphics[width=\linewidth]{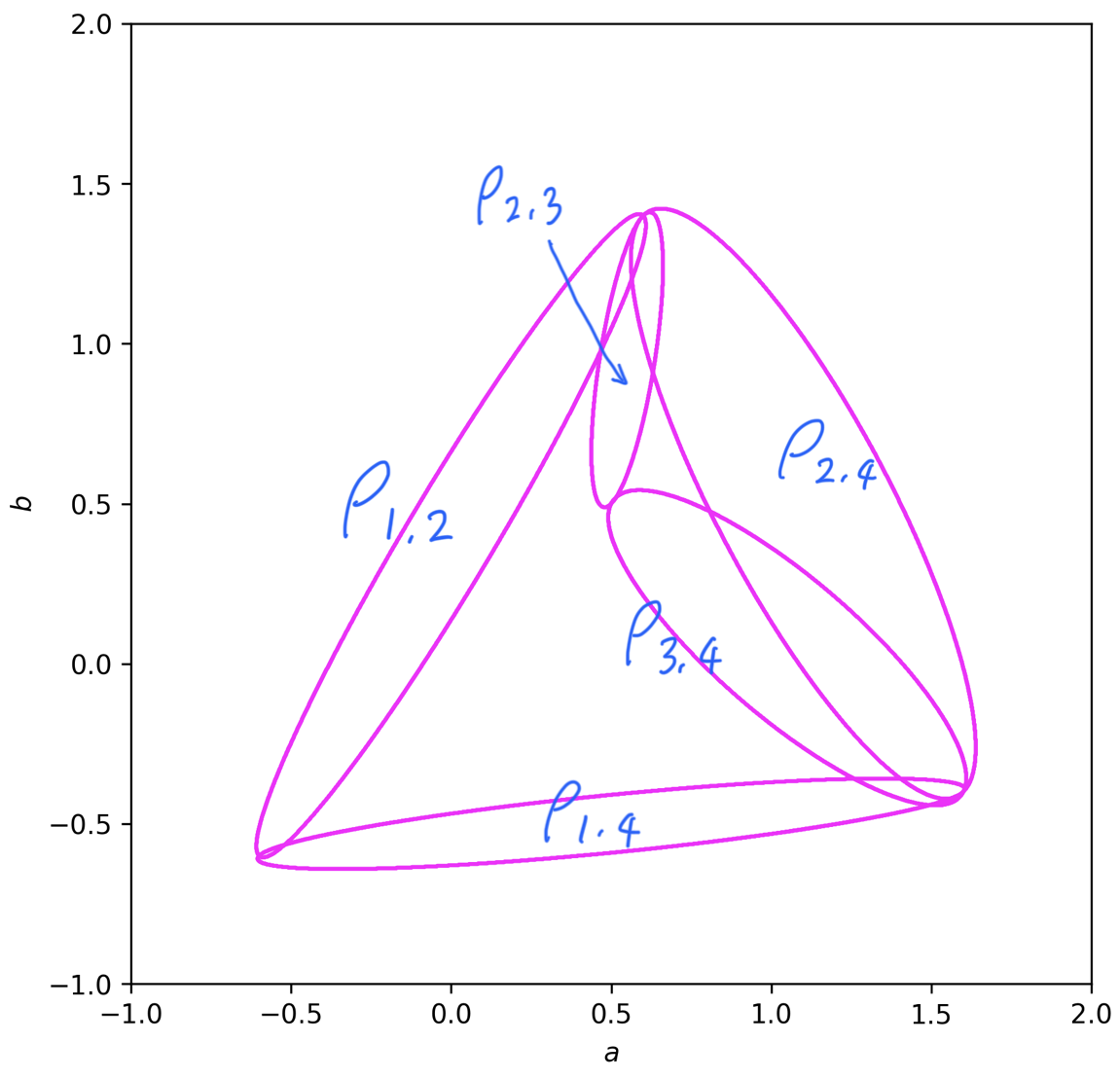}
        \caption{$\partial\,\tx{supp}(\rho_{approx})$ with each $\rho_{i,j}$ labeled.}
    \end{subfigure}
    \caption{The density $\rho$ and $\rho_{approx}$, error $\rho - \rho_{approx}$, and the boundary of $\tx{supp}(\rho_{approx})$ for Example~\ref{ex:sample_44}. }
    \label{fig:true_approx_error_boundary}
\end{figure}

\subsection{Support of Each Approximating Piece $\rho_{i,j}$}

Consider the approximating piece $\rho_{i,j}$ as defined in Theorem~\ref{thm:22_approx}. We shall prove the boundary of $\tx{supp}(\rho_{i,j}(\cdot,\cdot,\varepsilon_0))$ is an ellipse (for any fixed $\varepsilon_0$). We begin with some preliminary discussions. Define matrix $N$:
\begin{align*}
    N(a,b,\varepsilon, \lambda) := I - \frac{A(\varepsilon)}{a} - \lambda(bA(\varepsilon)-aB(\varepsilon)).
\end{align*}
We also consider the $2\times2$ principal submatrix $N_{\{i,j\}}$.

\begin{lemma}
\label{lemma:zeros_of_detN_superset_of_boundary}
For any $\varepsilon_0$, the zero set of $\tx{Disc}_{\lambda} \det N_{\{i,j\}}(\cdot, \cdot, \varepsilon_0, \lambda)$ equals to the zero set of $\tx{Disc}_{\lambda} \det(\lambda I - M^{\{i,j\}}(\cdot, \cdot, \varepsilon_0))$ and is a superset of $\partial \tx{supp}(\rho_{i,j}(\cdot,\cdot,\varepsilon_0))$
\end{lemma}
\begin{proof}
Fix any $\varepsilon_0$. Notice we have 
\begin{align*}
    \det(\lambda I - M^{\{i,j\}}(a, b, \varepsilon_0)) &:= \det\left(\lambda I - \left(I-\frac{A(\varepsilon_0)_{\{i,j\}}}{a}\right)(bA(\varepsilon_0)_{\{i,j\}}-aB(\varepsilon_0)_{\{i,j\}})^{-1}\right)\\
    &= \det\left(\lambda(bA(\varepsilon_0)_{\{i,j\}}-aB(\varepsilon_0)_{\{i,j\}}) - \left(I-\frac{A(\varepsilon_0)_{\{i,j\}}}{a}\right)\right) \det(bA(\varepsilon_0)_{\{i,j\}}-aB(\varepsilon_0)_{\{i,j\}})^{-1}\\
    &=: -\det N \det(bA(\varepsilon_0)_{\{i,j\}}-aB(\varepsilon_0)_{\{i,j\}})^{-1}.
\end{align*}
Hence for the discriminants, we have
\begin{align*}
    \tx{Disc}_{\lambda} \det(\lambda I - M^{\{i,j\}}(a, b, \varepsilon_0)) &= \tx{Disc}_{\lambda} \big(-\det N_{\{i,j\}} \det(bA(\varepsilon_0)_{\{i,j\}}-aB(\varepsilon_0)_{\{i,j\}})^{-1}\big)\\
    &= \det(bA(\varepsilon_0)_{\{i,j\}}-aB(\varepsilon_0)_{\{i,j\}})^{-2} \cdot \tx{Disc}_{\lambda} \det N_{\{i,j\}}.
\end{align*}
So an $(a,b)$ is a zero of $\tx{Disc}_{\lambda} \det(\lambda I - M^{\{i,j\}}(a, b, \varepsilon_0))$ if and only if it is a zero of $\tx{Disc}_{\lambda} \det N_{\{i,j\}}(a,b,\varepsilon_0, \lambda)$.\\

On the other hand, we know by Proposition~\ref{prop:zeros_of_Disc_contains_boundary} that the zeros of $\tx{Disc}_{\lambda} \det(\lambda I - M^{\{i,j\}}(\cdot, \cdot, \varepsilon_0))$ is a superset of $\partial\,\tx{supp}(\rho_{i,j}(\cdot,\cdot,\varepsilon_0))$.\\
\end{proof}

\newpage
\begin{proposition}
Consider any $i<j$ and $\varepsilon_0$. If $A(\varepsilon_0)$ and $B(\varepsilon_0)$ do not commute, then $\partial \tx{supp}(\rho_{i,j}(\cdot,\cdot,\varepsilon_0))$ is an ellipse with center $(\tx{tr} A/2, \tx{tr} B/2)$. Otherwise if $A(\varepsilon_0)$ and $B(\varepsilon_0)$ commute, then $\rho_{i,j}(\cdot, \cdot, \varepsilon_0)=0$ everywhere.
\end{proposition}
\begin{proof}
Fix $\varepsilon_0$ and assume $A(\varepsilon_0)$ and $B(\varepsilon_0)$ do not commute. We study $\tx{Disc}_{\lambda}\det N_{i,j}(\cdot, \cdot, \varepsilon_0, \cdot)$. Without loss of generality, assume $A(\varepsilon_0)$ is diagonal: otherwise we diagonalize $A = U^{-1}A_{diag}U$, where $U$ is orthogonal, and consider $A_{diag}, UBU^{-1}$ instead. The discriminant remains unchanged after this replacement:
\begin{align*}
    \det N &= \det(I - \frac{A}{a} - \lambda(bA-aB))\\
    &= \det(U^{-1})\det(I - \frac{A_{diag}}{a} - \lambda(bA_{diag} - aUBU^{-1}))\det U\\
    &= \det(I - \frac{A_{diag}}{a} - \lambda(bA_{diag} - aUBU^{-1})).
\end{align*}

Now we assume $A$ is diagonal. We have 
\begin{align*}
    N = \begin{bmatrix}
        1 - \frac{A_{1,1}}{a} - \lambda(bA_{1,1}-aB_{1,1}) & \lambda aB_{1,2}\\
        \lambda a B_{1,2} & 1 - \frac{A_{2,2}}{a} - \lambda(bA_{2,2}-aB_{2,2})
    \end{bmatrix}
\end{align*}
and
\begin{align*}
    &\det N = (1 - \frac{A_{1,1}}{a} - \lambda(bA_{1,1}-aB_{1,1}))(1 - \frac{A_{2,2}}{a} - \lambda(bA_{2,2}-aB_{2,2})) - \lambda^2 a^2 B_{1,2}^2\\
    &\tx{Disc}_{\lambda} \det N =: L(a,b)^2 + 4B_{1,2}^2(a - A_{1,1})(a-A_{2,2})\\
    &\qquad\qquad\quad = L(a,b)^2 + 4B_{1,2}^2\left(\left(a-\frac{A_{1,1}+A_{2,2}}{2}\right)^2 - \left(\frac{A_{2,2}-A_{1,1}}{2}\right)^2\right)
\end{align*}
\begin{align*}
    L(a,b) = -a(B_{2,2}-B_{1,1}) + b(A_{2,2}-A_{1,1}) + A_{1,1}B_{2,2}-A_{2,2}B_{1,1}.
\end{align*}
Hence,
\begin{align*}
    \tx{Disc}_{\lambda} \det N = 0 \iff \left(a-\frac{A_{1,1}+A_{2,2}}{2}\right)^2 + \left(\frac{L(a,b)}{2B_{1,2}}\right)^2 = \left(\frac{A_{2,2}-A_{1,1}}{2}\right)^2.
\end{align*}
Using the change of coordinate $\alpha = a - (A_{1,1}+A_{2,2})/2,\; \beta = L(a,b)/(2B_{1,2})$, we may write:
\begin{align*}
    &\tx{Disc}_{\lambda} \det N = 0 \iff \alpha^2 + \beta^2 = \left(\frac{A_{2,2}-A_{1,1}}{2}\right)^2\\
    &\tx{where } \begin{bmatrix}
        \alpha\\
        \beta
    \end{bmatrix} = 
    \begin{bmatrix}
        1 & 0\\
        -\frac{B_{2,2}-B_{1,1}}{2B_{1,2}} & \frac{A_{2,2}-A_{1,1}}{2B_{1,2}} 
    \end{bmatrix}
    \begin{bmatrix}
        a\\
        b
    \end{bmatrix} + 
    \begin{bmatrix}
        -\frac{A_{1,1}+A_{2,2}}{2}\\
        A_{1,1}B_{2,2}-B_{2,2}A_{1,1}
    \end{bmatrix}\\
    &\tx{with } \det \begin{bmatrix}
        1 & 0\\
        -\frac{B_{2,2}-B_{1,1}}{2B_{1,2}} & \frac{A_{2,2}-A_{1,1}}{2B_{1,2}} 
    \end{bmatrix} = \frac{A_{2,2}-A_{1,1}}{2B_{1,2}}.
\end{align*}
Recall our assumption that $A$ is diagonal and $AB-BA\neq 0$ implies that:
\begin{align*}
    AB - BA = \begin{bmatrix}
        0 & (A_{1,1}-A_{2,2})B_{1,2}\\
        (A_{2,2}-A_{1,1})B_{1,2} & 0 
    \end{bmatrix} \tx{ is not zero}
\end{align*}
hence the circle $\alpha^2 + \beta^2 = ((A_{2,2}-A_{1,1})/2)^2$ has non-zero radius. Consequently the zeros of $\tx{Disc}_{\lambda} \det N$ form an ellipse.\\

Lastly, we look for the center of the ellipse. We know the center is $(\alpha = 0, \beta = 0)$ after changing the coordinate. Equivalently, $(a-(A_{1,1}+A_{2,2})/2 = 0, L(a,b)/(2B_{1,2}) = 0)$. Solving for this gives $a = (A_{1,1}+A_{2,2})/2$ and $b = (B_{1,1}+B_{2,2})/2$. i.e., $a = \tx{tr}\,A/2, b = \tx{tr}\,B/2$.\\

Notice the formula for center aligns with our assumption that $A$ is diagonal: even if $A$ was not diagonal, after applying the diagonalization $A = UA_{diag}U^{-1}$, the trace is unchanged.
\end{proof}

\begin{remark}
Because each $\rho_{i,j}$ is a density as well, we know all statements that hold for a general $\rho$ holds for $\rho_{i,j}$ too, in particular the statements in the LSE setup.
\end{remark}

\subsection{Error in Terms of Density}

We prove the theorem claimed at the begining of Section~\ref{sec:22_approx_rho}:\\

\begin{figure}
    \centering
    \includegraphics[width=0.5\linewidth]{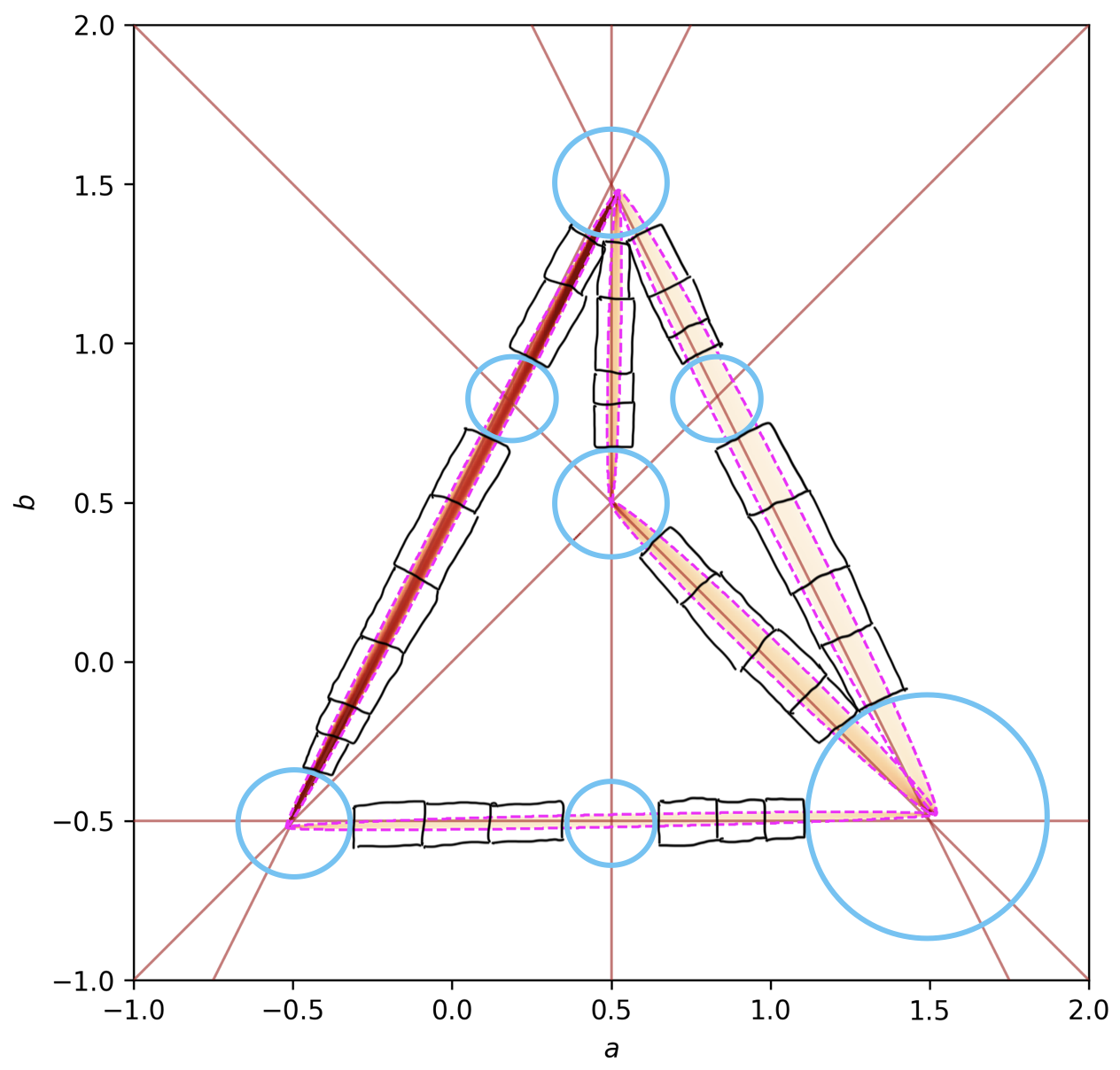}
    \caption{LSE covering of the support. The purple dashed lines are the boundary of $\tx{supp}(\rho_{approx}(\cdot, \cdot, \varepsilon))$ and the blue circles are the boundary of the open balls $B_{r_{i,j,k,l}}(L_{i,j}\cap L_{k,l})$. We cover the support of both $\rho_{approx}(\cdot, \cdot, \varepsilon)$ and $\rho(\cdot, \cdot, \varepsilon)$ by a finite number of the rectangles in the LSE setup.}
    \label{fig:LSE_covering}
\end{figure}

\textbf{Theorem~\ref{thm:22_approx}} [$2\times2$ Approximation of Density]\textbf{.} \textit{There exists $\varepsilon^*>0$ and $C>0$ such that for any $(a,b,\varepsilon) \in (\mathbb R^3 / \Omega_r) \times [0,\varepsilon^*]$, we have:}
\begin{align*}
    |\rho(a,b,\varepsilon) - \rho_{approx}(a,b,\varepsilon)| \leq C \varepsilon
\end{align*}
\textit{where the approximating density is defined as:}
\begin{align*}
    \rho_{approx} &:= \sum_{i<j} \rho_{i,j}\\
    &:= \sum_{i<j}\frac{1}{2\pi}\big(|\Im\, \lambda_1(M^{\{i,j\}}(a,b,\varepsilon))| + |\Im\, \lambda_2(M^{\{i,j\}}(a,b,\varepsilon))|\big)
\end{align*}
\textit{and} 
\begin{align*}
    M^{\{i,j\}}(a,b,\varepsilon) := \left(I-\frac{A_{\{i,j\}}(\varepsilon)}{a}\right)(bA_{\{i,j\}}(\varepsilon) - aB_{\{i,j\}}(\varepsilon))^{-1}
\end{align*}
\textit{where $A_{\{i,j\}}$ and $B_{\{i,j\}}$ refer to refer to the $2\times2$ principal submatrices of $A$ and $B$ corresponding to rows/columns $i$ and $j$, and $\Omega_r$ is defined in Definition~\ref{def:Omega}.}
\begin{proof}
First we prove the approximation locally in a LSE setup, then upgrade to a global statement by compactness of the support of $\rho$ and $\rho_{approx}$.\\

Consider the LSE setup and $\rho_{1,2}$. Notice $\rho_{1,2}$ and $\rho$ may be expressed in terms of the discriminants:
\begin{align*}
    &\rho_{1,2} = \frac{1}{2\pi} (|\Im \lambda_1| + |\Im \lambda_2|) = \frac{1}{\pi} |\Im\lambda_1|,\\
    &\tx{Disc}(\lambda I - M^{2\times2}) = (\lambda_1 - \lambda_2)^2 = -4|\Im\lambda_1|^2\\
    \implies& \rho_{1,2} = \frac{1}{2\pi} \sqrt{\max(-\tx{Disc}(\lambda I - M^{2\times2}), 0)}.
\end{align*}
Similarly for $\rho$:
\begin{align*}
    \rho = \frac{1}{2\pi}\sqrt{\max(-\tx{Disc}(\lambda I - T), 0)}.
\end{align*}
We also know for any $x, y \in \mathbb R$, $|\max(x,0)-\max(y,0)|\leq |x-y|$, and $|\sqrt{|x|} - \sqrt{|y|}|\leq \sqrt{|x-y|}$. Plug into density and discriminant:
\begin{align*}
    |\rho - \rho_{1,2}| &= \left|\sqrt{\max(-\tx{Disc}(\lambda I - T), 0)} - \sqrt{\max(-\tx{Disc}(\lambda I - M^{2\times2}),0)}\right|\\
    &\leq \left|\sqrt{-\tx{Disc}(\lambda I - T)} - \sqrt{-\tx{Disc}(\lambda I - M^{2\times2})} \right|\\
    &\leq \sqrt{|\tx{Disc}(\lambda I - T) - \tx{Disc}(\lambda I - M^{2\times2})|} \\
    &\leq C \varepsilon\sqrt{|a-1|+\varepsilon} \tag{Lemma~\ref{lemma:M_entry_Disc_bounds}(3)}\\
    &\leq C \varepsilon^{3/2}. \tag{Thm~\ref{thm:width_of_support_and_density}, Statement 1}
\end{align*}
We have finished the local part of the proof. Pick some $\varepsilon^*$ such that for every $\varepsilon \in [0,\varepsilon^*]$, we have:
\begin{enumerate}
    \item $\tx{supp}(\rho(\cdot, \cdot, \varepsilon)) \cup \tx{supp}(\rho_{approx}(\cdot, \cdot, \varepsilon))$ can be covered by a finite union of LSE rectangles, which is possible by the compactness of both supports.
    \item Outside the finite union of LSE rectangles, $\rho(\cdot, \cdot, \varepsilon) = \rho_{approx}(\cdot, \cdot, \varepsilon) = 0$, which is given by Theorem~\ref{thm:rho_locally_0_off_collision_line}.
\end{enumerate}
Figure~\ref{fig:LSE_covering} illustrates the intended covering. The claimed bound holds for this choice of $\varepsilon^*$. 
\end{proof}

The error of approximation for the $4\times4$ example, Example~\ref{ex:sample_44}, is plotted in the bottom-left panel of Figure~\ref{fig:true_approx_error_boundary}.

\subsection{Error in Terms of Boundary}
\label{sec:error_in_boundary}

We study the error between the boundary of $\rho$ and each $\rho_{i,j}$, say $\rho_{1,2}$. Again consider the LSE setup. We have concluded that the zeros of $\tx{Disc}(\lambda I - T)$ and $\tx{Disc}(\lambda I - M^{2\times2})$ is a superset of their respective boundaries, and in the non-degenerate case, the zeros equals to the boundaries (Statement (3) and (4) of Theorem~\ref{thm:width_of_support_and_density}). 

\begin{proposition}
\label{prop:error_of_boundary_non_degenerate}
(Error of Boundary, Non-Degenerate Case) Consider the LSE setup. There exists $C>0$ such that for any $b_0\in[-b^*, b^*]$ and any $\varepsilon\in [0,\varepsilon_{cvx}^*]$, if there are precisely two $a_- < a_+ \in \partial\,\tx{supp}(\rho(\cdot, b_0, \varepsilon))$ and two $a'_- < a'_+ \in \partial\,\tx{supp}(\rho_{1,2}(\cdot, b_0, \varepsilon))$, we are guaranteed that
\begin{align*}
    |a_- -a_+|\leq C \varepsilon^{3/2},\, |a'_- -a'_+|\leq C \varepsilon^{3/2}.
\end{align*}
\end{proposition}
\begin{proof}
Fix any $b_0\in [-b^*,b^*]$ and $\varepsilon_0\in[0,\varepsilon_{cvx}^*]$ and assume there are two distinct zeros for both boundaries. We know both $\tx{Disc}(\lambda I - T)(\cdot, b_0, \varepsilon_0)$ and $\tx{Disc}(\lambda I - M^{2\times2})(\cdot, b_0, \varepsilon_0)$ are strictly convex on $[1-a^*, 1+a^*]$ (Proposition~\ref{prop:Disc_at_most_two_zeros}). Say both second partials are lower bounded by some $K>0$. We also know that at their respective zeros, $\partial_v\tx{Disc}(\lambda I - T)(a_-, b_0, \varepsilon_0)<0$ and $\partial_v\tx{Disc}(\lambda I - T)(a_+, b_0, \varepsilon_0)>0$, similarly for $a_-'$, $a_+'$ and $\tx{Disc}(\lambda I - M^{2\times2})(\cdot, b_0, \varepsilon_0)$.\\

We prove the claimed bound on $a_+$ and $a_+'$, then similar proof holds for $a_-$ and $a_-'$. Fundamental theorem of Calculus gives that, for any $a$,
\begin{align*}
    &\partial_a\tx{Disc}(\lambda I - T)(a,b_0,\varepsilon_0) = \partial_a\tx{Disc}(\lambda I - T)(a_+,b_0,\varepsilon_0) + \int_{a_+}^a \partial_a^2 \tx{Disc}(\lambda I - T)(s,b_0,\varepsilon_0)\,ds,\\
    &\partial_a\tx{Disc}(\lambda I - M^{2\times2})(a,b_0, \varepsilon_0) = \partial_a\tx{Disc}(\lambda I - M^{2\times2})(a_+',b_0,\varepsilon_0) + \int_{a_+'}^a \partial_a^2 \tx{Disc}(\lambda I - M^{2\times2})(s,b_0,\varepsilon_0)\,ds.
\end{align*}
We use the first equation if $a_+ < a_+'$, and the second equation if $a_+> a_+'$; we skip the trivial case $a_+ = a_+'$. Say $a_+ < a_+'$. For any $a\in [a_+, a_+']$, we have:
\begin{align*}
    \partial_a\tx{Disc}(\lambda I - T)(a,b_0,\varepsilon_0) &\geq \partial_a\tx{Disc}(\lambda I - T)(a_+, \varepsilon_0) + K(a-a_+)\\
    &> K(a-a_+) > 0.
\end{align*}

We apply Fundamental Theorem of Calculus again on the discriminant:
\begin{equation}
\label{eq:Disc_lower_bounded_by_second_order}
\begin{split}
    \tx{Disc}(\lambda I - T)(a_+',b_0,\varepsilon_0) &= \tx{Disc}(\lambda I - T)(a_+,b_0,\varepsilon_0) + \int_{a_+}^{a_+'} \partial_v \tx{Disc}(\lambda I - T)(s,b_0,\varepsilon_0)\,ds\\
    &= 0 + \int_{a_+}^{a_+'} \partial_v \tx{Disc}(\lambda I - T)(s,b_0,\varepsilon_0)\,ds\\
    &\geq \frac{1}{2}K (a_+-a_+')^2 > 0.
\end{split}    
\end{equation}
We know the value of $\tx{Disc}(\lambda I - T)$ is controlled by Lemma~\ref{lemma:M_entry_Disc_bounds} (3):
\begin{align*}
    |\tx{Disc}(\lambda I - T)(a_+',b_0,\varepsilon_0)| = |\tx{Disc}(\lambda I - T)(a_+',b_0,\varepsilon_0) - \underbrace{\tx{Disc}(\lambda I - M^{2\times2})(a_+',b_0,\varepsilon_0)}_{=0}| \overset{Lemma~\ref{lemma:M_entry_Disc_bounds}(3)}{\leq} C\varepsilon^2(|a_+'|+\varepsilon_0).
\end{align*}
We also know by Proposition~\ref{prop:v_bound_by_eps_at_zero} that $|a|$ is bounded by first order of $\varepsilon$ at the zeros. i.e., there exists some constant, without loss of generality say the constant is $C$ as well, such that uniformly across all $\varepsilon_0\in [0,\varepsilon_{cvx}^*]$, we have
\begin{align}
\label{ineq:Disc_upper_bounded_by_third_order}
    |\tx{Disc}(\lambda I - T)(a_+',b_0,\varepsilon_0)| \leq C\varepsilon^2(|a_+'|+\varepsilon_0) \leq C\varepsilon_0^2 (1+C) \varepsilon_0 =: C_2 \varepsilon_0^3.
\end{align}
We have concluded:
\begin{align*}
    0 \overset{\eqref{eq:Disc_lower_bounded_by_second_order}}{\leq} \frac{1}{2}K (a_+'-a_+)^2 \overset{\eqref{eq:Disc_lower_bounded_by_second_order}}{\leq} \tx{Disc}(\lambda I - T)(a_+',b_0,\varepsilon_0) \overset{\eqref{ineq:Disc_upper_bounded_by_third_order}}{\leq} C_2 \varepsilon^3.
\end{align*}
Consequently $|a_+ - a_+'| \leq C_2 \varepsilon^{3/2}$.\\

Similar proof works for the other cases as well.\\
\end{proof}

\begin{figure}[H]
    \centering
    \includegraphics[width=0.5\linewidth]{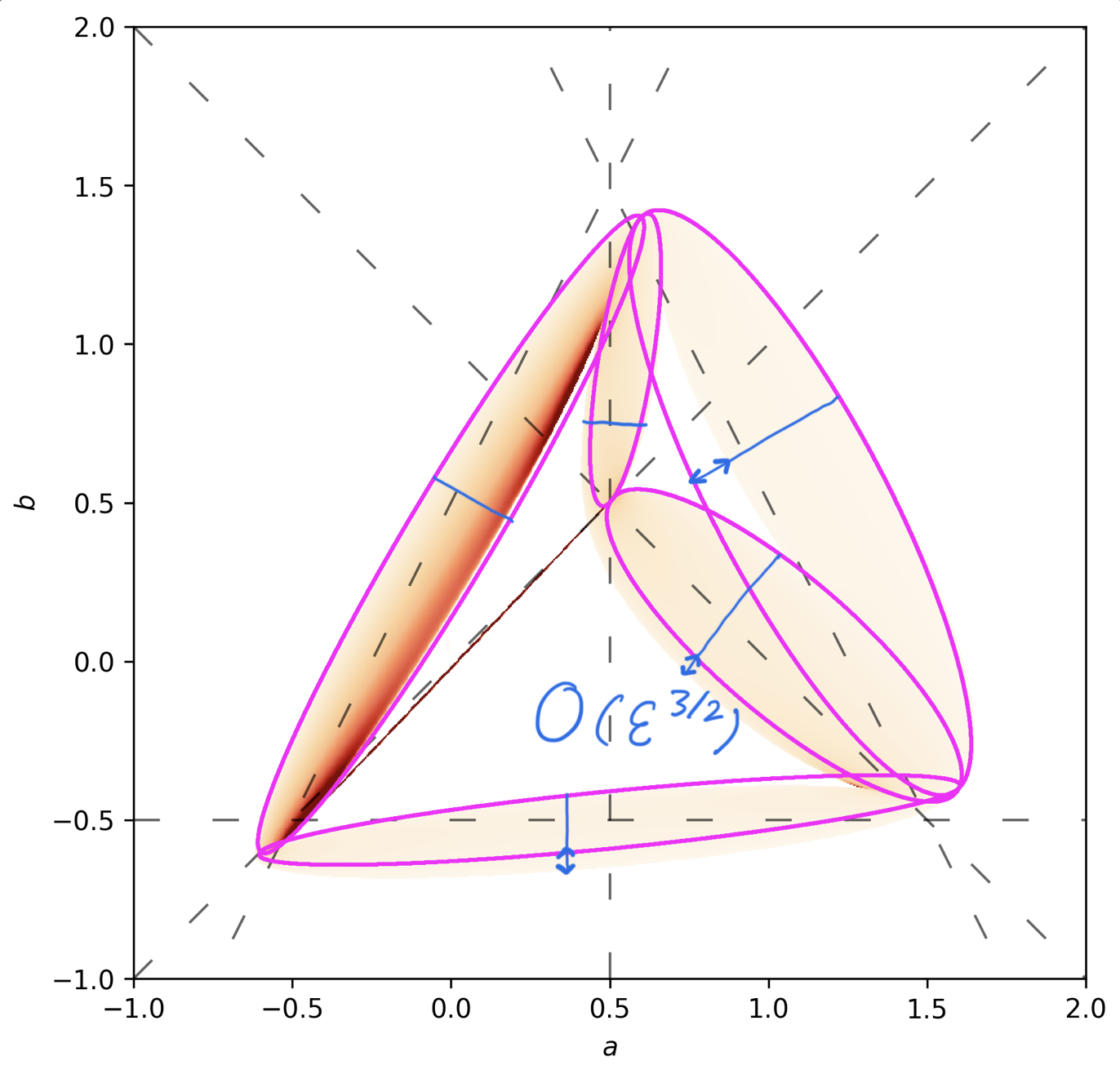}
    \caption{Error of boundary in the non-degenerate case. The distance between the boundary of $\tx{supp}(\rho)$ and $\tx{supp}(\rho)$, in the direction perpendicular to each collision line, is $O(\varepsilon^{3/2})$.}
\end{figure}

\begin{proposition}
(Error of Boundary, Degenerate Case) Consider the LSE setup. There exists $C>0$ such that the following holds for any $b_0\in[-b^*, b^*]$ and $\varepsilon_0\in[0,\varepsilon_{cvx}^*]$:
\begin{itemize}
    \item if $\rho(\cdot,b_0,\varepsilon_0) =0$ everywhere on $[1-a^*, 1+a^*]$, it must be the case that either $\rho_{1,2}(\cdot, b_0, \varepsilon_0) = 0$ everywhere as well, or the two $a_-'<a_+'\in \partial \tx{supp}(\rho_{1,2}(\cdot, b_0, \varepsilon_0))$ satisfy $|a_-' - a_+'| \leq C\varepsilon_0^{3/2}$.
    \item if $\rho_{1,2}(\cdot, a_0, \varepsilon_0) =0$ everywhere on $[1-a^*, 1+a^*]$, it must be the case that either $\rho(\cdot, b_0, \varepsilon_0) = 0$ everywhere as well, or the two $a_-<a_+\in \partial \tx{supp}(\rho(\cdot, b_0, \varepsilon_0))$ satisfy $|a_- - a_+| \leq C\varepsilon_0^{3/2}$.
\end{itemize}
\end{proposition}
\begin{proof}
We prove the second case of the first claim. The first case of the first claim is trivial; and the second claim holds by symmetry once the first one is proved.\\

Pick some $\varepsilon_0$ s.t the second case of the first claim holds. Say the minimum of $\tx{Disc}(\lambda I - M^{2\times2})$ is reached at $a_{min}'$. We know $a_-' < a_{min}' < a_+'$ by strict convexity. Again by strict convexity, we know $\partial_v \tx{Disc}(\lambda I - M^{2\times2})(\cdot, b_0, \varepsilon_0)$ is negative at $a_-'$, zero at $a_{min}'$, and positive at $a_+'$. Assume the second partial $\partial_a^2\tx{Disc}(\lambda I - M^{2\times2})(\cdot, b_0, \varepsilon_0)$ is lower bounded by some $K>0$. Then by Fundamental Theorem of Calculus, for any $a \in [a_{min}', a_+']$, we have:
\begin{align*}
    \partial_a \tx{Disc}(\lambda I - M^{2\times2})(a,b_0,\varepsilon_0) &\geq \partial_a \tx{Disc}(\lambda I - M^{2\times2})(a_{min}',b_0,\varepsilon_0) + K(a-a_{min}') \\
    &= K(a-a_{min}')> 0.
\end{align*}

By Fundamental Theorem of Calculus again:
\begin{align*}
    \tx{Disc}(\lambda I - M^{2\times2})(a_{min}',b_0,\varepsilon_0) &= \tx{Disc}(\lambda I - M^{2\times2})(a_+',b_0,\varepsilon_0) - \int_{a_{min}'}^{a_+'} \partial_a \tx{Disc}(\lambda I - M^{2\times2})(s,b_0,\varepsilon_0)\, ds\\
    &\leq -\frac{K}{2}(a_+' - a_{min}')^2 < 0.
\end{align*}

Lemma~\ref{lemma:M_entry_Disc_bounds}(4) gives a lower bound on $\tx{Disc}(\lambda I - M^{2\times2})(a_{min}',b_0,\varepsilon_0)$:
\begin{align*}
    -C \varepsilon^2(|a_{min}'| + \varepsilon_0) \leq \tx{Disc}(\lambda I - M^{2\times2})(a_{min}',b_0,\varepsilon_0) \leq -\frac{K}{2}(a_+' - a_{min}')^2 < 0,
\end{align*}
where the constant $C$ is uniform across all choices of $\varepsilon_0$ and $a_-', a_+'$. I.e.
\begin{align*}
    \sqrt{\frac{2C}{K}\varepsilon_0^2(|a_{min}'| + \varepsilon_0)} \geq (a_+' - a_{min}').
\end{align*}
Similarly we have $\sqrt{2C/K \varepsilon_0^2(|a_{min}'| + \varepsilon_0)} \geq (a_{min}' - a_-')$. We add up the two upper bounds and obtain:
\begin{align*}
    2\sqrt{\frac{2C}{K}\varepsilon_0^2(|a_{min}'| + \varepsilon_0)} \geq (a_+' - a_-').
\end{align*}
We know by Proposition~\ref{prop:v_bound_by_eps_at_zero} that both $|a_-'|$ and $|a_+'|$ are bounded by $C\varepsilon_0$, hence $|a_{min}'|$ is also bounded by  $C\varepsilon_0$:
\begin{align*}
    2\sqrt{\frac{4C}{K}\varepsilon_0^3} \geq (a_+' - a_-').
\end{align*}
So our claim holds. Similar proof works for the other case as well.
\end{proof}

\begin{figure}[H]
    \centering
    \includegraphics[width=0.48\linewidth]{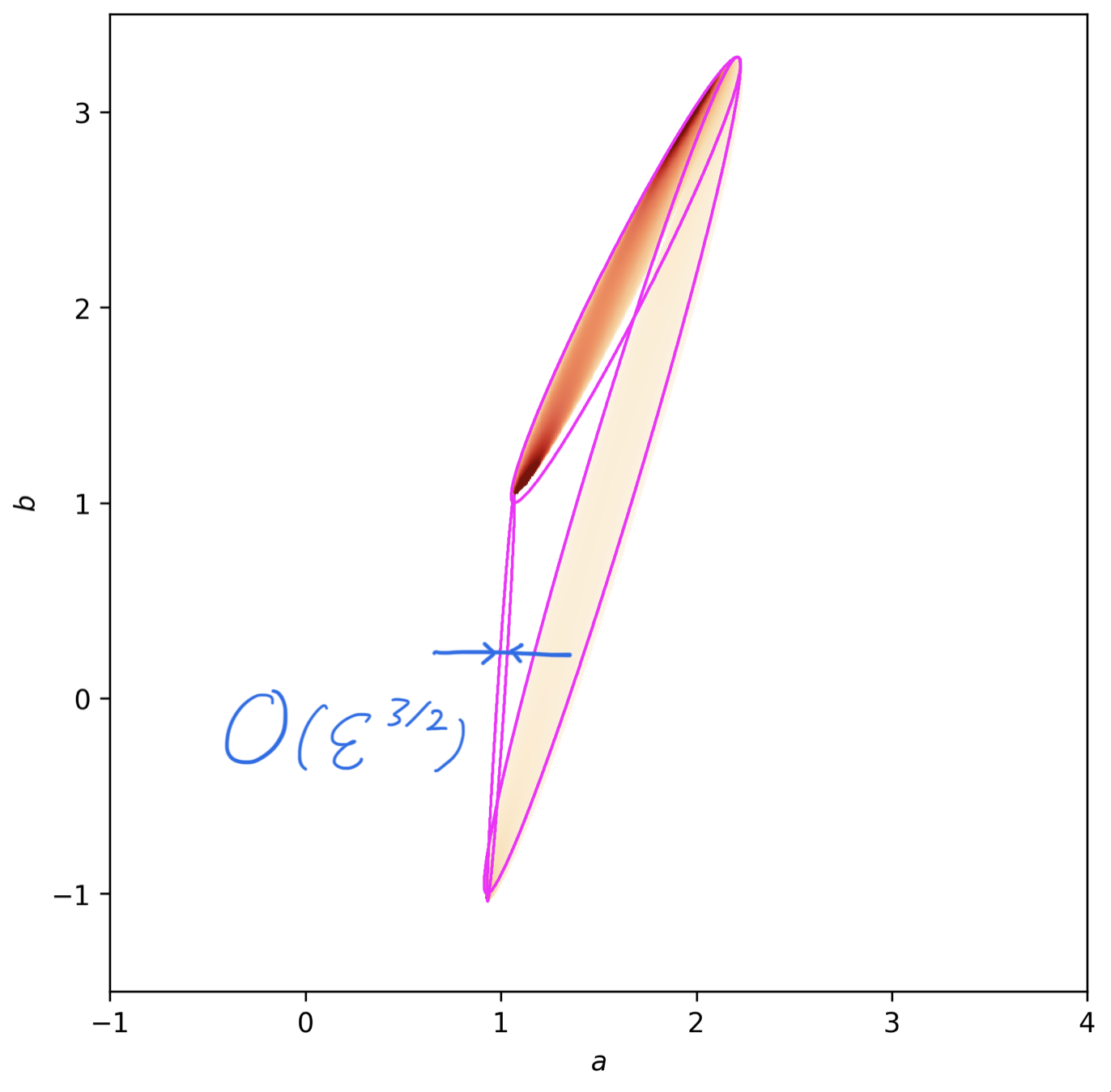}
    \hfill
    \includegraphics[width=0.48\linewidth]{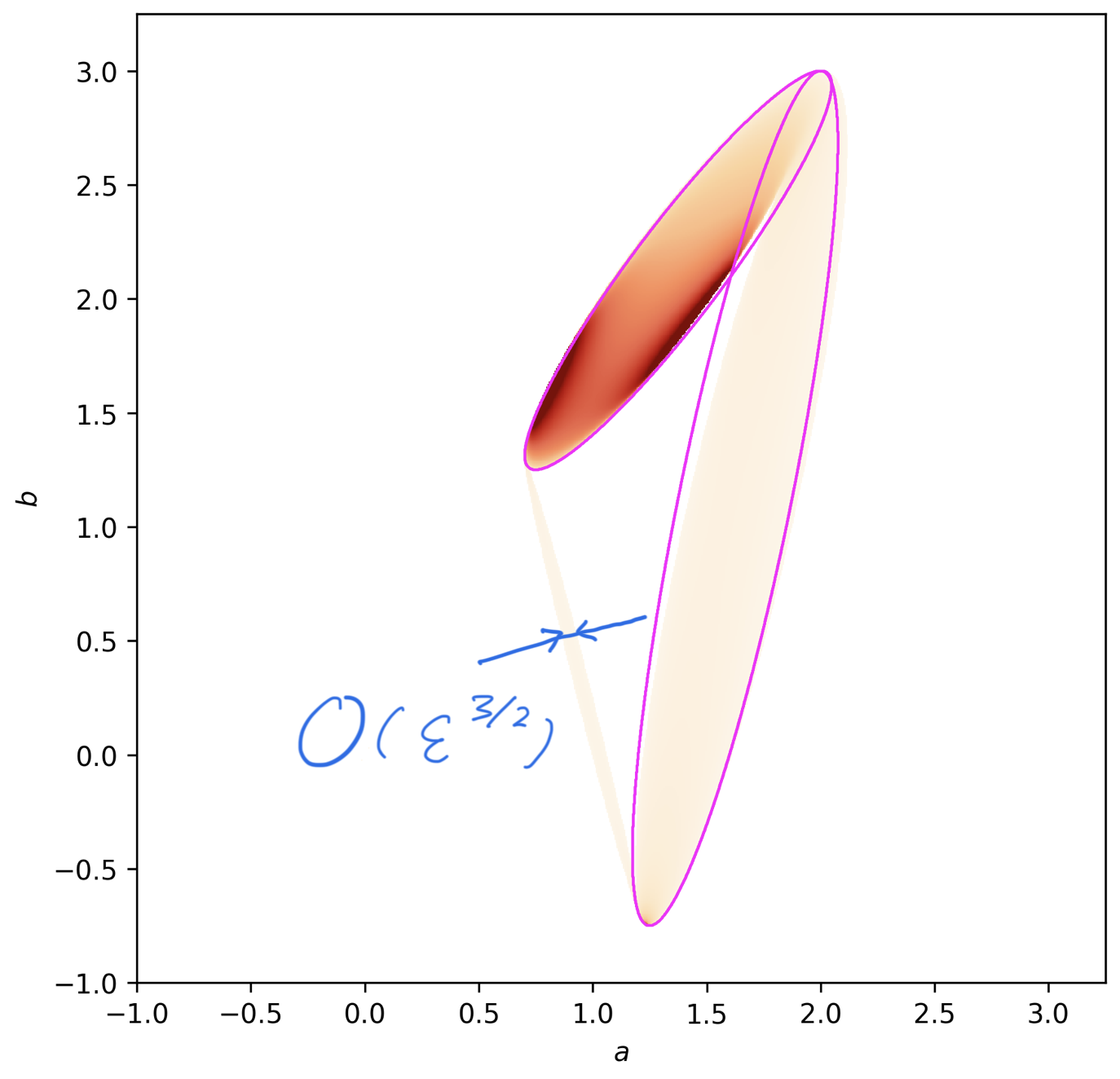}
    \caption{Error of boundary in the degenerate case. In the left panel, the true density $\rho$ degenerates and the width of $\tx{supp}(\rho_{approx})$ is $O(\varepsilon^{3/2})$; while we get the opposite in the right panel.}
\end{figure}

\section{Future Directions}
The present work leaves some natural questions unanswered. 
\begin{enumerate}
    \item Is the error bound $O(\varepsilon^{3/2})$ in Theorem~\ref{thm:22_approx} optimal? Our preliminary computation shows the error is likely second order.
    \item How to approximate $\rho$ near the intersection points of collision lines or in the absence of Assumption~\ref{asmp:non_collinear}, i.e., when the collision lines are not distinct?
    \item Whether projecting $(M(a,b,\varepsilon)$ to a larger spectral subspace can improve the approximation of $\rho$? In this work, we isolate the two eigenvalues near a collision line and construct a $2\times2$ approximation. Does isolating a cluster of $k>2$ eigenvalues provides a more precise approximation?
\end{enumerate}

\section*{Acknowledgements}

I want to thank my mentor Otte Hein{\"a}vaara. This research began as a Caltech SURF project and continued into the academic year. I want to thank my mentor for providing the idea of the project, the many instructions throughout the project, and revisions of this paper.\\

I also thank the Caltech Student-Faculty Programs Office and donors of the Ryser SURF Fellowship for sponsoring this project.\\

Lastly, several statements in this paper were proved with the help of ChatGPT. I appreciate the great convenience.

\printbibliography

@article{Hei_spectral,
  author  = {Hein{\"a}vaara, Otte},
  title   = {Tracial joint spectral measures},
  journal = {Inventiones mathematicae},
  year    = {2025},
  volume  = {239},
  number  = {2},
  pages   = {505--526},
  doi     = {10.1007/s00222-024-01309-6}
}

@phdthesis{Hei_spectral_thesis,
  author = {Hein{\"a}vaara, Otte},
  title  = {Tracial Joint Spectral Measures},
  school = {Princeton University},
  year   = {2024},
  month  = may,
  url    ={http://arks.princeton.edu/ark:/88435/dsp01z316q4969}
}

@book{HornJohnson2013,
  author    = {Horn, Roger A. and Johnson, Charles R.},
  title     = {Matrix Analysis},
  edition   = {2},
  publisher = {Cambridge University Press},
  year      = {2013}
}

@article{BauerFike1960,
  author  = {Bauer, F. L. and Fike, C. T.},
  title   = {Norms and Exclusion Theorems},
  journal = {Numerische Mathematik},
  volume  = {2},
  pages   = {137--141},
  year    = {1960},
  doi     = {10.1007/BF01386217}
}

@book{Barbeau1989,
  author    = {Barbeau, Edward J.},
  title     = {Polynomials},
  series    = {Problem Books in Mathematics},
  publisher = {Springer-Verlag},
  address   = {New York},
  year      = {1989}
}

@book{Kato1995,
  author    = {Kato, Tosio},
  title     = {Perturbation Theory for Linear Operators},
  edition   = {2},
  series    = {Classics in Mathematics},
  publisher = {Springer-Verlag},
  address   = {Berlin},
  year      = {1995},
  note      = {Reprint of the 1980 edition},
  doi       = {10.1007/978-3-642-66282-9}
}

\end{document}